\documentclass[11pt]{amsart}

\usepackage{fullpage}

\usepackage{url}

\usepackage{amssymb}

\usepackage{cite}

\renewcommand{\a}{\alpha}
\renewcommand{\b}{\beta}
\newcommand{\g}{\gamma}
\renewcommand{\d}{\delta}
\newcommand{\D}{\Delta}
\newcommand{\e}{\varepsilon}
\newcommand{\f}{\varphi}
\newcommand{\s}{\sigma}
\renewcommand{\S}{\Sigma}
\renewcommand{\k}{\kappa}
\renewcommand{\l}{\lambda}

\renewcommand{\O}{\Omega}
\renewcommand{\o}{\omega}
\newcommand{\G}{\Gamma}

\newcommand{\cO}{{\mathcal O}}
\newcommand{\cF}{{\mathcal F}}
\newcommand{\cC}{{\mathcal C}}
\newcommand{\cM}{{\mathcal M}}

\newcommand{\cS}{{\mathcal S}}
\newcommand{\cB}{{\mathcal B}}

\newcommand{\cE}{{\mathcal E}}
\newcommand{\cU}{{\mathcal U}}
\newcommand{\cV}{{\mathcal V}}
\newcommand{\cP}{{\mathcal P}}
\newcommand{\cN}{{\mathcal N}}

\newcommand{\cD}{{\mathcal D}}

\newcommand{\cK}{{\mathcal K}}

\newcommand{\cZ}{{\mathcal Z}}
\newcommand{\cA}{{\mathcal A}}

\newcommand{\bC}{\mathbb C}
\newcommand{\bR}{\mathbb R}

\newcommand{\bZ}{\mathbb Z}

\newcommand{\bH}{\mathbb H}

\newcommand{\bP}{\mathbb P}

\newcommand{\Di}{\mathrm{Diff}}

\newcommand{\Co}{\mathrm{Conf}(S^{2})}

\newcommand{\be}{\begin{equation}}
\newcommand{\ee}{\end{equation}}

\renewcommand{\to}{\rightarrow}

\renewcommand{\phi}{\varphi}
\renewcommand{\epsilon}{\varepsilon}
\renewcommand{\hat}{\widehat}

\newcommand{\<}{\langle}
\renewcommand{\>}{\rangle}

\newcommand{\w}{\widetilde}

\theoremstyle{plain}
\newtheorem{theorem}{Theorem}[section]

\newtheorem{remark}[theorem]{Remark}

\newtheorem{lemma}[theorem]{Lemma}

\newtheorem{proposition}[theorem]{Proposition}

\theoremstyle{definition}

\def\endproof{\qed \medskip}
\def\blacksquare{\hbox to .60em {\vrule width .60em height .60em}}

\numberwithin{equation}{section}

\begin{document}

\title[ ]{The Nirenberg problem of prescribed Gauss curvature on $S^{2}$}

\author[ ]{Michael T. Anderson}

\begin{abstract}
We introduce a new perspective on the classical Nirenberg problem of understanding the possible Gauss curvatures 
of metrics on $S^{2}$ conformal to the round metric. A key tool is to employ the smooth Cheeger-Gromov compactness 
theorem to obtain general and essentially sharp {\it a priori} estimates for Gauss curvatures $K$ contained in naturally 
defined stable regions. We prove that in such stable regions, the map $u \to K_{g}$, $g = e^{2u}g_{+1}$ is a proper Fredholm 
map with well-defined degree on each component. This leads to a number of new existence and non-existence results. 
We also present a new proof and generalization of the Moser theorem on Gauss curvatures of even conformal metrics 
on $S^{2}$. 

   In contrast to previous work, the work here does not use any of the Sobolev-type inequalities of 
Trudinger-Moser-Aubin-Onofri. 
\end{abstract}

\address{Department of Mathematics, Stony Brook University, Stony Brook, N.Y.~11794-3651, USA} 
\email{anderson@math.sunysb.edu}
\urladdr{http://www.math.sunysb.edu/$\sim$anderson}

\thanks{Partially supported by NSF grant DMS 1607479}

\maketitle

\setcounter{section}{0}
\setcounter{equation}{0}

\section{Introduction}

In this work, we study the well-known Nirenberg problem: which smooth functions $K$ on $S^{2}$ are realized as the Gauss curvature 
of a metric $g$ on $S^{2}$ pointwise conformal to the standard round metric $g_{+1}$ of radius 1 on $S^{2} \subset \bR^{3}$? 
For $g = e^{2u}g_{+1}$, the equation for the Gauss curvature $K$ of $g$ is 
\be \label{u}
e^{2u}K = 1 - \D u.
\ee
so that the Nirenberg problem asks to characterize for which $K$ is the nonlinear PDE \eqref{u} solvable. It is also of interest to 
understand the uniqueness or multiplicity of solutions. 

   Prior work on the Nirenberg problem is based on its variational formulation initiated by Moser \cite{M1}. 
For any given $K$ consider the functional 
\be \label{J}
J_{K}(u) = \int_{S^{2}} (|d u|^{2} + 2u)dv_{0} - \log(\int_{S^{2}} K e^{2u}dv_{0}), 
\ee
where $dv_{0}$ is the usual round area form on $S^{2}$, normalized to unit area. 
Critical points $u$ of $J$ are given by solutions of \eqref{u}, (up to constants), i.e.~\eqref{u} is the Euler-Lagrange equation of 
$J_{K}$. Moser in \cite{M2} began the analysis of $J_{K}$ based on fundamental Sobolev inequalities due to Trudinger \cite{T} 
and Moser \cite{M1}, leading then to a solution of the existence problem for even functions $K$ on $S^{2}$. 

  The functional $J_{K}$ does not satisfy the Palais-Smale Condition C and except in case $K = const$, the infimum of $J_{K}$ is 
never realized, \cite{Ho}. Most all general existence results after \cite{M1} have been obtained by identifying conditions on 
$K$ that prevent the {\it a priori} possible blow-up or bubble behavior of suitable approximate solutions of \eqref{u}. 
These approximate solutions are sequences approaching a minimax (mountain pass) critical point of the functional $J_{K}$. 
Further Sobolev inequalities of Onofri \cite{O} and Aubin \cite{Au} play an important role in this analysis. 

  Fundamental progress on the Nirenberg problem was made by Kazdan-Warner \cite{KW}, \cite{KW1}, Aubin \cite{Au}, 
Chang-Yang \cite{CY1}-\cite{CY3}, Chang-Gursky-Yang \cite{CGY}, Chen-Ding \cite{CD} with further important progress 
made by Han \cite{ZH}, Chang-Liu \cite{CL} and many others. The full literature on this problem is vast and we mention 
here only further work in \cite{Ho}, \cite{Ji}, \cite{St}; see also \cite{JLX} for a recent survey and further references.  
The problem is still far from a complete resolution, even in the case of positive Gauss curvature. 

  In this paper, we give new proofs of almost all of the results above from a different perspective, and derive as well new 
existence and non-existence results. In particular, the proofs of these results do not rely on any of the sharp Sobolev 
inequalities or (with one exception) any variational formulation of the problem. 

 To introduce the point of view taken here, we recast the problem \eqref{u} as follows. Let 
$C^{m,\a}$, $m \geq 4$, $\a \in (0,1)$ be the Banach space of $C^{m,\a}$ functions $u: S^{2} \to \bR$. This is 
considered as the space of conformal factors for metrics $g = e^{2u}g_{+1}$. Consider the map 
$$\pi: C^{m,\a} \to C_{+}^{m-2,\a},$$
\be \label{pi}
\pi(u) = K_{g},
\ee
where $K_{g}$ is the Gauss curvature of $g$ and $C_{+}^{m-2,\a}$ is the space of 
$C^{m-2,\a}$ functions $K$ such that $K(x) > 0$ for some $x \in S^{2}$. By integrating \eqref{u}, one easily sees 
that $\pi$ maps into $C_{+}^{m-2,\a}$. The Nirenberg problem then asks to describe the image $Im \, \pi$ of $\pi$. 
It is well-known that \eqref{pi} is not solvable for all $K \in C_{+}^{m-2,\a}$. There is a fundamental obstruction due to 
Kazdan-Warner \cite{KW}: if $K$ is the Gauss curvature of the conformal metric $g$, then 
\be \label{kw}
\int_{S^{2}}X(K)dv_{g} = 0,
\ee
for all linear vector fields $X$ on $S^{2}$; $X = \nabla \ell$ where $\ell$ is the restriction of a linear function on 
$\bR^{3}$ to $S^{2}$. For example, the functions $K = 1 + \ell$ are not realizable as the Gauss curvature of a 
conformal metric $g = e^{2u}g_{+1}$.

   From \eqref{u} and standard elliptic theory, $\pi$ is a smooth (nonlinear) Fredholm map, of Fredholm index 0. 
A broad understanding of the global properties of the map $\pi$, for instance its image, requires understanding when 
the Fredholm map $\pi$ is proper; recall that a map is proper if the inverse image of any compact set in the target is 
compact in the domain. Loosely speaking, $\pi: \cU \to \cV$ proper on a domain $\cV$ is equivalent to the statement 
that control of $K = \pi(u) \in \cV$ implies control of any solution $u$ of \eqref{u}; this equivalent to the existence 
of {\it a priori} estimates for solutions of \eqref{u} with $K \in \cV$. 

\medskip 

   Now a general principle in Riemannian geometry is that control of the curvature $K$ implies control of the metric $g$, 
at least given suitable bounds on global quantities such as volume and diameter. This principle, which has origins in early 
work in global differential geometry, is expressed concretely in two fundamental results in geometry, namely  the Cheeger 
finiteness theorem \cite{Ch} and the (smooth) Gromov convergence theorem \cite{Gr}, \cite{GW}, \cite{Pe}, now generally 
referred to as the Cheeger-Gromov convergence theorem; this is discussed in detail in Section 2. The control on the metric 
given by the control on the curvature comes however only {\it modulo diffeomorphisms}; the action of diffeomorphisms 
is crucial and cannot be avoided. 

  Applied to the case of $S^{2}$, it follows that control of $K$ implies control of $g = e^{2u}g_{+1}$ modulo diffeomorphisms 
of $S^{2}$, given suitable control on the area and diameter of $(S^{2}, g)$. Now the Nirenberg problem per se is of course 
{\it not} diffeomorphism invariant. However, on $S^{2}$ there is a unique conformal class, and so all metrics may be 
pulled back by diffeomorphisms to the fixed conformal class $[g_{+1}]$. This pull-back is unique modulo the conformal 
group $\Co$.  

  This highlights the central role played by the non-compact conformal group $\Co$ of $S^{2}$. The importance of the 
action of $\Co$ is certainly well-known from all previous analyses of the Nirenberg problem from the purely PDE point of view. 
However, the emphasis is quite different here, taking the Riemannian geometry of the metrics much more closely into 
consideration. 

   To describe the first main result, we consider first the behavior of $K$ near its zero level set. Thus, suppose $q \in S^2$ 
satisfies $K(q) = 0$ and that in a small neighborhood of $q$, $K$ has the form 
\be \label{criterion}
K = \chi r^{2a},
\ee
where $r(x) = dist_{S^2(1)}(x, q)$, $a \geq 1$ and $\chi$ is at least continuous near $q$. Thus we are considering points on $\{K = 0\}$ 
which are either non-degenerate local maxima or minima of $K$ (when $a = 1$ and $\chi(q) \neq 0$) or arbitrary degenerate critical points 
(when $a > 1$ or $\chi(q) = 0$). Let 
\be \label{P}
\cC = \{K \in C_{+}^{m-2,\a}: {\rm if} \ \eqref{criterion} \ {\rm holds \ at} \ q, \ {\rm then} \ \D K(q) > 0\}.
\ee
In particular, if $K \in \cC$ and $K(q) = 0$, then $q$ is not a local maximum point of $K$. We note that if $q$ is a 
non-degenerate saddle point of $K$, then \eqref{criterion} does not hold near $q$; such $K$ are thus in $\cC$. 
We also note that it is straightforward to see that $\cC$ is invariant under the 
action of the diffeomorphism group $\Di^{m-2,\a}$ on $C_{+}^{m-2,\a}$. Next let  
\be \label{cN2}
\cN = \{K \in C_{+}^{m,\a}: |\nabla K|(p) + |\D K|(p) > 0, \forall p \ {\rm s.t.} \ K(p) > 0\}.
\ee
Thus $K \in \cN$ if at any critical point $q$ of $K$ with $K(q) > 0$, $\D K(q) \neq 0$. When $K > 0$ everywhere, 
$\cN$ is exactly the space of non-degenerate functions as defined by Chang-Yang in \cite{CY3}, \cite{CGY}. 

   In contrast to $\cC$, the space $\cN$ is {\it not} invariant under the action of the diffeomorphism group 
$\Di^{m-2,\a}$. It is easy to see that both $\cN$ and $\cC$ are open and dense in $C_{+}^{m-2,\a}$ and are 
invariant under the action of $\Co$. Let 
$$\cU = \pi^{-1}(\cC \cap \cN),$$
and, setting 
$$\cK = \cC \cap \cN,$$ 
let 
\be \label{pi0}
\pi_{0} = \pi|_{\cU}: \cU \to \cK, \ \ \pi_{0}(u) = K_{u}, 
\ee
be the restriction of $\pi$ to $\cU$. The first main result is the following: 
  
\begin{theorem}
On the domain $\cU$, the curvature map 
$$\pi_{0}: \cU \to \cK$$
in \eqref{pi0} is a proper Fredholm map of index 0. 
\end{theorem} 

  In particular, for any $K \in \cK$, the space of solutions $\pi^{-1}(K)$ of \eqref{u} is compact and is finite for generic 
$K \in \cK$. Theorem 1.1 is also essentially sharp in that $\pi$ does not extend to a proper map on any larger open domain 
$\cV$ containing $\cU$, cf.~Remark 5.11. Namely, we show that bubble solutions of unbounded area may form as 
$K \in \cK$ approaches $\partial \cC$. This follows from the basic work of Ding-Liu in \cite{DL} and later by Borer-Galimberti-Struwe 
in \cite{BGS},  cf.~Remark 2.7. Similarly, we show that solutions in $\cU$ degenerate due to the non-compactness of the conformal 
group $\Co$ as $K \in \cK$ approaches $\partial \cN$. This in turn follows from the basic work of Chang-Yang \cite{CY1}-\cite{CY3}. 

  Theorem 1.1 gives {\it a priori} estimates for solutions $u$ of \eqref{u} for all $K = \pi(u) \in \cK$, cf.~Remark 2.10. 
Such estimates have been previously obtained by Chang-Gursky-Yang \cite{CGY} for $K \in \cK$ in case $K > 0$ and 
extended by Chen-Li in \cite{CL} to a larger, but still strictly proper, subset of $K \in \cK$ with $K$ of variable sign. 
Theorem 1.1 gives the first such estimates in the full region $\cK$. Moreover, the proof of Theorem 1.1 is quite different 
than previous work on this topic. As noted above, it does not rely on any delicate analysis related to sharp Sobolev inequalities 
of Trudinger-Moser-Aubin type or on any blow-up analysis or proof by contradiction. The proof is arguably much simpler 
than previous proofs. 

  We prove in Section 2 that the complement $\partial \cK = C_{+}^{m-2,\a} \setminus \cK$ is a closed, rectifiable set of 
codimension one in $C_{+}^{m-2,\a}$. Generically $\partial \cK$ is a collection of smooth, codimension one hypersurfaces.

\medskip 

  It is well-known that proper Fredholm maps $F: X \to Y$ of index zero between Banach spaces $X$,$Y$, (or orientable 
Banach manifolds with $Y$ connected) have a well-defined $\bZ$-valued degree given by 
$$deg \, F = \sum_{x \in F^{-1}(y)}sign(DF_{x}) \in \bZ$$
wbere $y$ is any regular value of $F$ and the sign of $DF_{x}$ is $\pm 1$ according to whether $DF(x)$ preserves or reverses 
orientation at $x$. 

  Thus let $\cK^{i}$ denote the collection of path components of $\cK \subset C_{+}^{m-2,\a}$ and let $\cU^{i} = 
(\pi_{0})^{-1}(\cK^{i})$ be the inverse image in $C^{m,\a}$. Note that {\it a priori} $\cU^{i}$ may not be connected. 
The restriction of $\pi_{0}$ to $\cU^{i}$ gives a proper Fredholm map  
\be \label{pii}
\pi^{i}: \cU^{i} \to \cK^{i},
\ee
with a well-defined degree. To calculate the degree, let $K$ be any regular value of $\pi^{i}$ in $\cK^{i}$ and without 
loss of generality, assume $K$ is a Morse function on $S^{2}$. Let 
\be \label{zi}
Z_{i}(K) = \{q:\nabla K(q) = 0, \D K(q) < 0 \ {\rm and} \ K(q) > 0 \},
\ee
so that $Z_{i}(K)$ is a finite set of points in $S^{2}$. Let $ind \,q$ be the index of the non-degenerate critical point 
$q$ of $K$; ($ind \,q$ is even for $q$ a local minimum or maximum, while $ind \,q$ is odd if $q$ is a saddle point). 
The degree of $\pi^{i}$ on the domain of positive curvature functions $K > 0$ in $\cK^{i}$ was calculated by 
Chang-Gursky-Yang \cite{CGY} and we prove that the same formula extends to the full domains $\cK^{i}$. 

\begin{theorem} 
One has 
\be \label{degform0}
deg \,\pi^{i} = \sum_{q \in Z_{i}(K)} (-1)^{ind \, q} - 1,
\ee
for any $Z_{i}(K)$ as in \eqref{zi}. Equivalently, if $M$ is the number of positive local maxima $K$ and $s^{-}$ is the 
number of saddle points of $K$ where $K > 0$ and $\D K < 0$, then 
$$deg \,\pi^{i} = M-s^{-}-1.$$
\end{theorem}

  It is not difficult to see that $deg \, \pi^{i}$ may assume any value in $\bZ$, i.e.~for each $n \in \bZ$, there 
exists $i$ such that $deg \, \pi^{i} = n$, cf.~Lemma 4.1.  It is not clear (and is perhaps unlikely) whether the 
components $\cK^{i}$ of $\cK$ are uniquely determined by their degree, cf.~Remark 4.2 for further discussion. 

  Combining the variational formulation \eqref{J} of equation \eqref{u} with Theorem 1.1 leads to a second formula for the degree. 

\begin{theorem}
   The degree of $\pi^{i}$ is also given (up to a fixed sign determined by orientation) by 
\be \label{deg2}
deg \,\pi^{i} = \sum_{u_{j} \in (\pi^{i})^{-1}(K)}(-1)^{ind(u_{j})},
\ee
where $K$ is any regular value of $\pi^{i}$.
\end{theorem}
Here the index $ind(u_{j})$ is the dimension of the negative eigenspace of the Hessian $D^{2}J_{K}$ at $u_{j}$, 
for $J_{K}$ as in \eqref{J}.

\medskip 

  Theorems 1.2 or 1.3 prove existence of solutions of \eqref{u}, together with a signed count on the number of 
generic solutions, for all $K \in C_{+}^{m-2,\a}$ except $K$ in the degree zero components and $K$ on the 
boundary $K \in \partial \cK$. We prove in Remark 4.2 that the degree $0$ component is connected and so call this 
component $\cK^{0}$. Of course all $K$ satisfying \eqref{kw} are necessarily in $\cK^{0}$ or $\partial \cK$. 

  To obtain further information in these regions, we analyse the structure of the regular and singular points of 
the Fredholm map $\pi$. The Sard-Smale theorem \cite{Sm} implies that the regular values of $\pi$ are open and dense 
in the target space $C_{+}^{m-2,\a}$. Of course, {\it a priori}, on the degree zero component 
\be \label{p0}
\pi^{0}: \cU^{0} \to \cK^{0},
\ee
there may be no corresponding regular points of $\pi$; a generic value in the target space $\cK^{0}$ may have 
empty inverse image. (For instance, the image may be of high codimension or even empty). 

\begin{theorem} The set $\S$ of singular points of $\pi$ in \eqref{pi} is a stratified space with strata of codimension $s \geq 1$ in 
$C^{m,\a}$. In particular, the regular points and regular values of $\pi$ are open and dense in the domain and range of $\pi$. 
\end{theorem}

We refer to Section 5 for the definition of a stratified space. This result gives a sharpening of the Perturbation Theorem 
in \cite{KW1}. The image of $\pi^{0}$ in \eqref{p0}, 
$$\cE = Im \, \pi^{0} \subset \cK^{0} \subset C_{+}^{m-2,\a}$$ 
is a non-empty {\it closed} subset of $\cK^{0}$. This is a consequence of the fact that $\pi^{0}$ is proper, 
cf.~Proposition 5.7 below. The Kazdan-Warner obstruction \eqref{kw} implies that $\pi^{0}$ is not surjective, so the 
complement of the image $\cD = (Im \, \pi^{0})^{c} \subset \cK^{0}$ is a non-empty open set in $\cK^{0}$ 
(and hence open in $C_{+}^{m-2,\a}$). We show in Proposition 4.4 that for any given $K$, if $\ell$ is any linear 
function with $|\ell|$ sufficiently large (depending on $K$), then \eqref{u} is not solvable for $K + \ell$, i.e.~$K 
+ \ell \notin Im \, \pi$. 

In the region $\cE$, solutions to the Nirenberg problem typically come in pairs, with 
$J_{K}$-index of opposite parity. We prove in Proposition 5.8 that the boundary $\cB \subset \partial \cE$ separating 
the existence region $\cE$ from the non-existence region $\cD$ in $C_{+}^{m-2,\a}$ is generically a smooth bifurcation 
locus for the creation or annihilation of pairs of solutions of the Nirenberg problem, 
cf.~also Remark 5.9.  

   Next, the boundary $\partial \cK$ decomposes as a union 
$$\partial \cK =  \partial \cC \cup \partial \cN.$$
The first region corresponds to a wall-crossing of a local maximum of $K$ passing through the value $0$ while 
the second region corresponds to a wall-crossing where the value of $\D K$ changes sign at a critical point 
$q$ of $K$ with $K(q) > 0$. As noted above, degenerations given by bubbles of infinite area may form on 
approach to $\partial \cC$ while the non-compactness of the conformal group $\Co$ causes degeneration on 
approach to $\partial \cN$.

Regarding $\partial \cC$, we have the following result: 
\begin{theorem} 
For any two distinct components $\cK^{m}$, $\cK^{n}$ of degree $m$ and $n$ respectively with $m \neq 0$ and $m < n$, 
all Morse functions $K \in \cM$ in $\partial \cC_{m}^{n} = \partial \cC \cap \partial \cK^{m}\cap \partial \cK^{n}$ are 
realizable as Gauss curvatures of conformal metrics, i.e.
$$\partial \cC_{m}^{n} \cap \cM \subset Im \, \pi.$$ 
\end{theorem}
This is proved in Proposition 4.7 below. A result of this generality is not known for $\partial \cN$. 
We prove a partial result holds near the point $K = 1$, cf.~ the discussion preceding Remark 5.11. Also large 
families of explicit solutions in $\partial \cN$ are constructed in Proposition 4.5 and Remark 4.6; we show that 
if $u$ is any eigenfunction of the Laplacian on $S^{2}(1)$ then $\pi(u) = K_{u} \in \partial \cN$.

\medskip 

  We conclude the paper with an existence result for functions $K$ which have a symmetry which breaks the action of 
the conformal group $\Co$ on $C_{+}^{m-2,\a}$, as in Moser's theorem \cite{M2} on the existence of solutions of \eqref{u} 
with $K$ even. Thus, let $\G$ be a finite subgroup of $O(3) = Isom(S^{2}(1))$. The action of $\G$ is said to break the 
(non-compact) action of the conformal group if the $\G$-action and $\Co/O(3)$-action on $S^{2}$ do not commute. Thus 
for any conformal dilation $\f \in \Co/O(3) \simeq \bR^{3}$, with $\f \neq Id$, there exists $\g \in \G$ such that 
$$\g\circ \f \neq \f \circ \g.$$
For instance, it is easy to see that the $\bZ_{2}$-action generated by the antipodal map breaks the non-compact 
action of $\Co$. 

\begin{theorem} 
Let $\G \subset O(3)$ be a finite subgroup breaking the action of $\Co$. Then any $\G$-invariant $C_{+}^{m-2,\a}$ 
function $K: S^{2} \to \bR$ is the Gauss curvature of a $\G$-invariant conformal metric $g$ on $S^{2}$. 
\end{theorem} 

\smallskip 

  The contents of the paper are briefly as follows. In Section 2, we prove Theorem 1.1 using the geometric convergence 
theory of Cheeger-Gromov, combined with a basic estimate of Astala \cite{As} on regularity of quasi-conformal mappings 
and the basic work of Chang-Yang \cite{CY3}. We also relate the decomposition \eqref{pii} into path 
components of $\cK$ to the more well-known decomposition of the space of Morse functions in $C_{+}^{m-2,\a}$ into 
Morse chambers. Theorems 1.2 and 1.3 are proved in Section 3, building on the basic degree formula of Chang-Gursky-Yang 
\cite{CGY} for Theorem 1.2 and on general degree theoretic properties for proper Fredholm maps for Theorem 1.3. 
Section 4 is a bridge between the earlier and later sections and discusses several existence and non-existence results, 
including Theorem 1.5, in the regions $\cK^{0}$, $\partial \cC$ and $\partial \cN$. In Section 5, we analyse the structure 
of singular and regular points of $\pi$ and prove Theorem 1.4.   Finally, Section 6 is devoted to the proof of Theorem 1.6.  
Further results on the topics above are given in the individual sections. 

\medskip 

  I am most grateful to Alice Chang, and also to Paul Yang and Matt Gursky, for discussions which greatly helped clarify  
their previous work on this topic to me and which led to significant corrections of an earlier version of this paper. I am 
also very grateful to the referee for the very detailed comments and careful work on the paper. 

\section{Properness of $\pi$.} 

  The basic issue in studying the global properties of the Fredholm map $\pi$ in \eqref{pi} is whether $\pi$ is proper. 
Thus, given any sequence $g_{i} = e^{2u_{i}}g_{+1}$ of metrics in the standard round conformal class with 
Gauss curvature $K_{i} = K_{g_{i}}$ such that 
\be \label{prop}
K_{i} \to K \ \ {\rm in} \ \ C_{+}^{m-2,\a},
\ee
the issue is when does a subsequence of $\{u_{i}\}$ converge to a limit in $C^{m,\a}$. It is well-known that this is 
not true in general, due to the non-compactness of the conformal group $\Co \simeq PGL(2, \bC)$ of fractional 
linear transformations of $(S^{2}(1), g_{+1})$. The maximal compact subgroup of $PGL(2, \bC)$ is $O(3)$, the isometry 
group of $S^{2}(1)$, and the quotient satisfies 
\be \label{h3}
PGL(2, \bC)/O(3) \simeq \bH^{3} \simeq \bR^{3}.
\ee
Elements in the quotient $\bR^{3}$ are identified with the conformal dilations as follows. Given a point $p 
\in S^{2}$, in the chart given by stereographic projection from the pole $p$,  
$$\f_{p,t}(x) = tx,$$
represents the conformal dilation by $t$ with source $-p$ and sink $p$. One has $\f_{p,1} = Id$, which serves as the 
origin $0 \in \bR^{3}$. Also $(\f_{p,t})^{-1} = \f_{p, t^{-1}} = \f_{-p, t}$. Similarly the map $\f_{p,q,t} = 
t(x - x_0)$ represents the conformal dilation by $t$ with source $q$ and sink $p$, where $q$ maps to $x_0$ under 
stereographic projection. (This corresponds to a conjugation of $\f_{p,t}$ by a rotation). 

  The action (or skew-action) of $\Co$ on the target space $C_{+}^{m-2,\a}$ of curvature functions is by pre-composition, 
$$(K, \f) \to K \circ \f.$$
This action is clearly not proper, since the constant functions $K = const$ are fixed points of the action. 
On the other hand, it is easy to see that the action of $\Co$ is proper on $C_{+}^{m-2,\a}$ away from the 
constant functions. 

  To describe the action of $\Co$ on the domain space $C^{m,\a}$ of conformal factors, let 
$\f \in {\rm Conf}(S^{2})$ be an orientation preserving conformal diffeomorphism, so that $\f^{*}g_{+1} = 
\chi_{\f}^{2}g_{+1}$, where $\chi_{\f}^{2} = det D\f > 0$. For $u \in C^{m,\a}$ with $g = e^{2u}g_{+1}$, 
one has 
$\f^{*}g = e^{2u\circ \f}\f^{*}g_{+1} = e^{2u\circ \f}\chi_{\f}^{2}g_{+1}$. 
Thus ${\rm Conf}(S^{2})$ acts on $C^{m,\a}$ as 
\be \label{confact}
(\f, u) \to \f^{*}u = u_{\f} = (u \circ \f) + \log \chi_{\f}.
\ee
For conformal dilations $\f_{p,t}$, the volume distortion $\chi_{p,t}^{2}$ tends to the Dirac measure supported 
at $p$ as $t \to \infty$. Note that since $\f_{-p,s}\circ \f_{p,t} = \f_{p,\frac{t}{s}}$, $\f_{-p,s}^*(\f_{p,t}^*(0)) = \log \chi_{p, \frac{t}{s}}$ so that 
\be \label{compf}
\f_{-p,s}^*(\log \chi_{p,t}) = \log \chi_{p,t}\circ \f_{-p,s} + \log \chi_{-p,s} = \log \chi_{p, \frac{t}{s}}.
\ee
For later purposes (e.g.~Lemma 2.5), observe that the level set $\{\log \chi_{p,t} = 0\}$ converges to the point $p$ as $t \to \infty$ while 
the level set $\{\f_{-p,s}^*\log \chi_{p,t} = 0\}$ remains bounded away from $p$ and $-p$ as $t, s \to \infty$ but $\frac{t}{s}$ remains 
bounded. 

  In contrast to the action of $\Co$ on the target space $C_{+}^{m-2,\a}$ of curvature functions, the action of $\Co$ 
on the domain $C^{m,\a}$ is smooth and proper; all orbits are properly embedded. The map $\pi$ is equivariant with respect 
to these actions of $\Co$, i.e. 
\be \label{comm1}
\pi(\f^{*}u) = \f^{*}\pi(u).
\ee
In particular, $\pi$ maps $\Co$ orbits in the domain $C^{m,\a}$ to $\Co$ orbits in the target $C_{+}^{m-2,\a}$ and so 
$\pi$ is not proper onto a neighborhood of the constant functions $K = const$. 

Following \cite{Au}, let $\cS$ be the space of conformal factors $e^{2u}$ with zero center of mass, 
\be \label{cS}
\cS = \{ u \in C^{m,\a}: \forall i, \int_{S^2} e^{2u}x_{i} dv_{+1} = 0\},
\ee
where $x_{i}$ are the restrictions of the standard linear coordinate functions on $\bR^{3}$ to $S^{2}(1)$. The space 
$\cS$ is a smooth codimension 3 hypersurface in $C^{m,\a}$. It is proved in \cite{CL} that $\cS$ is a smooth global 
slice for the non-compact action of $\Co$ on $C^{m,\a}$, in that there is a smooth diffeomorphism 
\be \label{cslice}
C^{m,\a} \simeq \cS \times (\Co / O(3)) \simeq \cS \times \bR^{3}.
\ee
In particular, for any $u \in C^{m,\a}$ there is a unique conformal dilation $\f_{p,t}$, $(p,t) \in \bR^{3}$, such that 
$(\f_{p,t})^{*}u \in \cS$. The dilation $\f_{p,t}$ depends smoothly on $u$. 

\medskip 

  Given this background, we study the properness issue in general by first recalling the Cheeger-Gromov convergence theorem, 
expressed here in the very special case of $S^{2}$. Let $W^{k,p}$ denote the Sobolev space of functions on $S^{2}$ with $k$ 
weak derivatives in $L^{p}$; convergence in the weak topology in $W^{k,p}$ is denoted by $W_*^{k,p}$.  

\begin{proposition} {\rm (Cheeger-Gromov)}. 
Let $g_{i}$ be a sequence of $C^2$ metrics on $S^{2}$ with uniformly bounded Gauss curvature $K_{i} = K_{g_{i}}$, 
\be \label{kbound}
|K_{i}| \leq K_{0} < \infty.
\ee
 Suppose there are constants $v_{0} > 0$ and $D_{0} < \infty$ such that 
\be \label{ad}
area_{g_{i}}S^{2} \geq v_{0} \ \ {\rm and} \ \ diam_{g_{i}}S^{2} \leq D_{0}.
\ee
There there is a sequence of diffeomorphisms $\psi_{i} \in \Di^{2,\a}\cap W^{3,p}$ such that 
\be \label{dpsi}
\w g_{i} = \psi_{i}^{*}g_{i}
\ee
is uniformly bounded in $C^{1,\a}\cap W^{2,p}$, for any $\a < 1$, $p < \infty$. In particular, a subsequence converges 
\be \label{wconv}
\w g_{i} \to \w g \ \ {\rm in} \ \ C^{1,\a},
\ee
and weakly in $W^{2,p}$ (i.e.~in $W_*^{2,p}$) for any $p < \infty$ to a limit $C^{1,\a}\cap W^{2,p}$ metric $\w g$ on $S^{2}$. 
\end{proposition}

{\bf Proof:} This is a standard result that holds in all dimensions and detailed proofs are given in \cite{GW}, \cite{Pe}, \cite{Pet} for 
instance. For later purposes, we sketch some of the main ideas of the proof. A key idea is to locally represent any metric 
$g$ by well-controlled harmonic coordinate charts, i.e.~local diffeomorphisms 
$$F: B(r_{0}) \subset \bR^{2} \to S^{2},$$
with $F^{-1} = (x_{1}, x_{2})$ given by a pair of harmonic functions $x_{i}$ with respect to $g$ (or equivalently with 
respect to $g_{+1}$ when $g$ is conformal to $g_{+1}$). The local charts $F$ are constructed so that the metric 
$\bar g = F^{*}g$ is bounded in the $C^{1,\a}\cap W^{2,p}$ norm by the $L^{\infty}$ norm of $\bar K = K \circ F$, given the 
global bounds \eqref{ad}. The key to this control is the well-known expression for the components $\bar g_{\a\b}$ 
of a metric in local harmonic coordinates given by 
\be \label{harm}
\D_{\bar g}\bar g_{\a\b} + Q_{\a\b}(\bar g, \partial \bar g) = -2(Ric_{\bar g})_{\a\b} = -2\bar K\bar g_{\a\b},
\ee
where $Q$ is a lower order term, quadratic in $\bar g$ and its first derivatives. Here $\D_{\bar g}$ is the Laplacian of 
$\bar g$, given in the harmonic coordinate chart as 
$$\D_{\bar g} = \bar g^{ab}\partial_{a}\partial_{b}.$$
The system \eqref{harm} is an elliptic system for $\bar g_{\a\b}$ and the $C^{1,\a}\cap W^{2,p}$ bounds on $\bar g_{\a\b}$ 
follow from standard elliptic estimates, given $L^{\infty}$ control on $\bar K$. 

The diffeomorphisms $\psi_{i}$ in \eqref{dpsi} are constructed by assembling the overlap maps of such local 
charts by a geometric gluing, such as a center of mass averaging. Elliptic regularity again shows that each $\psi_{i}$ is 
$C^{2,\a}\cap W^{3,p}$ smooth. The uniform $C^{1,\a}\cap W^{2,p}$ local control on the metrics $F_{i}^{*}g_{i}$ then 
easily leads to the convergence \eqref{wconv} by the Arzela-Ascoli theorem. 

{\endproof} 

  We may apply Proposition 2.1 to a sequence $g_i = e^{2u_i}g_{+1}$ of conformal metrics on $S^2$ satisfying 
\eqref{kbound}-\eqref{ad}. The uniformization theorem on $S^{2}$ implies that any $C^{1,\a}\cap W^{2,p}$ metric $\w g$ is of the form 
$\w g = \psi^{*}(e^{2u}g_{+1})$ for some $\psi \in \Di^{2,\a}\cap W^{3,p}$ and $u \in C^{1,\a}\cap W^{2,p}$. 
Thus, in a subsequence, 
\be \label{a11}
\w g_{i} = \psi_{i}^{*}(g_{i}) \to \w g = \psi^{*}(e^{2u}g_{+1}),
\ee
in $C^{1,\a} \cap W_{*}^{2,p}$. Replacing  $\psi_{i}$ by $\w \psi_{i} = \psi_{i}\circ \psi^{-1}$ and dropping the tilde from the notation, 
\eqref{a11} then becomes 
\be \label{a31}
\psi_{i}^{*}(e^{2u_{i}}g_{+1}) \to e^{2u}g_{+1},
\ee
in $C^{1,\a} \cap W_{*}^{2,p}$. The equation \eqref{a31} means that the diffeomorphisms $\psi_{i}$ are 
$(1 + \e_{i})$-quasi-conformal with $\e_{i} \to 0$ as $i \to \infty$. By a standard result (cf.~\cite{LV}) any 
sequence of $(1+\e_{i})$-quasi-conformal diffeomorphisms has a subsequence converging to a limit, 
modulo the action of the conformal group. (This is the normal family principle for quasi-conformal mappings). 
Any such limit is clearly a conformal map of $S^{2}$. Thus, there exist conformal maps $\f_{i} \in \Co$ such that 
\be \label{id1}
\eta_{i} : = \f_{i}^{-1}\circ \psi_{i} \to Id.
\ee
The diffeomorphisms $\eta_{i}$ are $(1+\e_{i})$ quasi-conformal with $\e_{i} \to 0$ as $i \to \infty$. The equations 
\eqref{a31}-\eqref{id1} will be a main focus of attention after Remark 2.7 below. First however we must address the 
situation when one of the hypotheses in \eqref{ad} fails.

\begin{remark}
{\rm  We first note that Proposition 2.1 is essentially local; suitable versions of the result hold for bounded domains or 
complete, non-compact manifolds, when one works with respect to given base points, cf.~\cite{Pet} for further discussion 
of the general situation. 

In the case at hand, consider a pointed sequence $(S^2, g_i, y_i)$, $g_i = e^{2u_i}g_{+1}$ satisfying \eqref{kbound}, where $y_i$ 
is any sequence of base points $y_i \in S^2$. Suppose the following local non-collapse assumption (analogous to \eqref{ad}) holds: 
\be \label{nc}
area_{g_i}B_{y_i}(1) \geq v_0 > 0.
\ee
Here $v_0$ is an arbitrary positive constant and $B_{y_i}(1)$ denotes the $g_i$-geodesic ball of radius $1$ about $y_i$. 
Then the Cheeger-Gromov theory implies that a subsequence converges in the pointed $C^{1,\a}\cap W_*^{2,p}$ topology modulo 
diffeomorphisms to a complete Riemannian surface $(\O, \w g, \w y)$. In general $\O$ is an open surface; if $\O$ is closed 
(and hence compact), then necessarily $\O = S^2$. If $\O$ is open, then $diam_{g_i}S^2 \to \infty$. The convergence 
modulo diffeomorphisms means that there is an exhaustion $\O_k$ of $\O$, with $\O_k \subset \O_{k+1}$, $\cup \O_k = \O$, 
$\w y \in \O_k$ for all $k$, and embeddings $\psi_k: \O_k \to S^2$, $\psi_k(\w y) = y_k$, such that $\w g_i = \psi_k^*(e^{2u_i}g_{+1}) 
\to \w g|_{\O_k}$. In addition, by the uniformization theorem, there is a conformal embedding $F: \O \to S^2$ 
of the abstract space $(\O, \w g)$ so that $\w g = F^*(e^{2u}g_{+1})$ for some function $u$ defined on ${\rm Im}(F)$. The mappings 
$\psi_k$ and $F$ are not unique; one may compose them with diffeomorphisms for instance. 

  Passing to a diagonal subsequence of $\{i, k\}$ and relabeling, it follows from the above that 
\be \label{oconv}
(F^{-1})^*\psi_i^*(e^{2u_i}g_{+1}) \to e^{2u}g_{+1},
\ee
in $C^{1,\a}\cap W_*^{2,p}$ uniformly on compact subsets of $F(\O)$, analogous to \eqref{a31}. As above, the maps $\psi_i\circ F^{-1} $ form a sequence of 
quasi-conformal diffeomorphisms with dilatation $1+\e_i \to 1$ and as in \eqref{id1}, there exists a sequence of conformal 
transformations $\f_i$ such that 
\be \label{toId}
\f_i \circ \psi_i \circ F^{-1} \to Id,
\ee
on ${\rm Im}(F)$. 

  This shows that an analog of the discussion following Proposition 2.1 holds when the hypothesis \eqref{ad} is replaced by the 
weaker non-collapse assumption \eqref{nc}. 

  Finally, we note that if one assumes $g_i \in C^{m+2}$ and a stronger bound on the covariant derivatives of the Gauss curvatures, 
\be \label{covbound}
|\nabla^{j}K_{i}|_{g_{i}} \leq K_{m},
\ee 
$1 \leq j \leq m$, then the convergence in \eqref{wconv} or \eqref{oconv} above can be improved to convergence in $C^{m+1,\a} \cap  
W_*^{m+2,p}$. Note the bound \eqref{covbound} is invariant under diffeomorphism. On the other hand, the convergence 
$K_{i} \to K$ in $C^{m-2,\a}$ as in \eqref{prop} is not invariant under diffeomorphism. Such convergence does not imply 
convergence of $K_{i}\circ \psi_{i}$ to a limit $\w K$ in, say, $C^{\a}$, if the diffeomorphisms $\psi_{i}$ tend to infinity. 

}
\end{remark}

  The Gauss-Bonnet theorem gives 
$$4\pi = \int_{S^{2}}K_{g}dv_{g} \leq \max |K_{g}| area_{g}S^{2},$$
so that the lower area bound in \eqref{ad} (but not necessarily \eqref{nc}) holds automatically under uniform curvature bounds. 
If one has a uniform positive lower bound 
\be \label{k0}
K_{g} \geq K_{0} > 0
\ee
on the Gauss curvature, then the diameter bound in \eqref{ad} is also automatic by the well-known Bonnet-Myers theorem. 
In general, when \eqref{k0} does not hold for some $K_0 > 0$, there may be sequences for which $diam_{g_{i}}S^{2} \to \infty$. 

\medskip 

   We show next that the diameter must remain uniformly bounded when the metrics $g$ are in a 
fixed conformal class with curvature $K$ in a compact subset of $\cC$, for $\cC$ defined as in 
\eqref{P}. Let $[g]$ denote the pointwise conformal class of a metric $g$. 

\begin{theorem}  Let $g_{i}$ be a sequence of $C^{m,\a}$ metrics on $S^2$ with $m \geq 4$. If $[g_{i}] = [g_{+1}]$ and 
$K_{i} \in \cC$ with  
\be \label{maxK}
K_{i} \to K \in \cC \subset C_{+}^{m-2,\a}, 
\ee
then there is a constant $D_{0}$ such that 
\be \label{diamb}
diam_{g_{i}}S^{2} \leq D_{0} < \infty,
\ee
and hence $area_{g_{i}}S^{2} \leq A_{0}$, for some $A_{0} < \infty$. 
\end{theorem}

   The proof of this result is rather long and broken into the following lemmas and propositions. The first proposition does not require 
the limit curvature $K \in \cC$. 

\begin{proposition}
Suppose $[g_i] = [g_{+1}]$ with $K_i \to K$ in $C^{m-2,\a}$ and suppose in addition there is a constant $A_{0} < \infty$ such that 
\be \label{abound}
area_{g_{i}}S^{2} \leq A_{0}.
\ee
Then \eqref{diamb} holds. 
\end{proposition}

{\bf Proof:} The proof is by contradiction and so we assume there is a sequence $g_{i}$ as above satisfying \eqref{abound} 
but for which 
\be \label{dinf}
diam_{g_{i}}S^{2} \to \infty. 
\ee
For any $\e > 0$ and any metric $g$ on $S^{2}$ one has a thick-thin decomposition of 
$(S^{2}, g)$, 
\be \label{tt}
S^{2} = U^{\e} \cup U_{\e},
\ee
where $U^{\e}$ is the set of points $x$ where the injectivity radius $inj_{g}(x) \geq \e$ while $U_{\e}$ is the complement. 
We first claim that if \eqref{abound} and \eqref{dinf} hold, then both sets $U^{\e}$ and $U_{\e}$ are non-empty, for any given $\e > 0$ 
sufficiently small, provided $i$ is sufficiently large. 

  To see this, we recall the collapse theory of Cheeger-Fukaya-Gromov \cite{CFG} which states that if a manifold $(M, g)$ satisfies 
$|Rm| \leq \Lambda$ and $inj_{g}M \leq \e_{0}$, then $M$ has an $\cF$-structure determined by the collection of short geodesic 
loops based at any point $p \in M$. Here $\e_{0}$ is a fixed constant, depending only on $\Lambda$ and the dimension $n$ of $M$ 
while $Rm$ is the Riemann curvature of $(M, g)$.  
In two dimensions, $\cF$-structures are particularly simple; for any $p \in M$, there is a unique short geodesic loop based at $p$, 
of length equal to $2inj_{g}(p)$, (assuming $inj_{g}(p)$ is sufficiently small compared with $\Lambda^{-1}$). Such loops may be 
smoothed at the base points to give a foliation of $M$ by circles. In particular any surface with an $\cF$-structure has Euler 
characteristic $\chi(M) = 0$; (the same holds in all dimensions by \cite{CFG}). Since $\chi(S^{2}) \neq 0$, it follows that 
$$U^{\e} \neq \emptyset,$$
for all $\e \leq \e_{0}$, where $\e_{0}$ depends only on $\max |K|$. On the other hand, if $U_{\e} = \emptyset$, then every point of 
$(S^{2}, g_{i})$ has injectivity radius at least $\e$. Since the curvature $K$ is uniformly bounded, standard comparison geometry, 
cf.~\cite{Pet} for example, shows that $area_{g_{i}}(B_{p}(1)) \geq a_{0}$, where $a_{0}$ depends only on $\e$ and $\max |K|$. The area bound 
\eqref{abound} then implies a uniform upper bound on the number of disjoint geodesic 1-balls $B_{p_{j}}(1)$ which can be contained in 
$(S^{2}, g_{i})$. Taking a maximal collection of such balls, the balls $B_{p_{j}}(2)$ of radius 2 then cover the manifold $(S^{2}, g_{i})$ 
which implies that $diam_{g_{i}}S^{2}$ is uniformly bounded above. This contradicts \eqref{dinf}. This proves the claim above. 

  Let $\w U^{\e}$ be the unit tubular neighborhood of $U^{\e}$, $\w U^{\e} = \{x \in (S^2, g): dist(x, U^{\e}) \leq 1\}$. Again, standard 
comparison geometry, cf.~\cite{Gr} or \cite{Pet}, implies that $\w U^{\e} \subset U^{\l\e}$, for a (small) constant $\l$ depending only on a bound 
for $\max |K|$. As above, since each connected component of $\w U^{\e}$ has a uniform lower bound on its area, it follows that $\w U^{\e}$ 
has a bounded number of components (independent of $i$). By passing to a subsequence, we may assume $\w U^{\e}$ has a 
fixed number of components (for any given choice of $\e \leq \e_{0}$). Further, the same covering argument as above shows 
there is a uniform upper bound on the diameter of each component of $\w U^{\e}$; 
$$diam_{g_{i}}\w U^{\e} \leq D,$$
where $D$ depends only on $\e$ and $\max |K_{i}|$.  

  Let $\w U_{\e} = S^2 \setminus \w U^{\e}$, so that $\w U_{\e} \subset U_{\e}$. In the following, we assume $\e$ is chosen so that 
$\e \leq \e_{0}$, so that each component of $\w U_{\e}$ has an $\cF$-structure; more precisely each component of $\w U_{\e}$ 
has a small thickening $\w U_{\e}' \supset \w U_{\e}$ which has an $\cF$-structure so that $\w U_{\e}'$ has a foliation by disjoint short 
circles. We relabel $\w U_{\e}'$ to $\w U_{\e}$ and similarly relabel $\w U^{\e}$ to $S^2 \setminus \w U_{\e}$. Thus each component 
of $\w U_{\e}$ is topologically an annulus $I\times S^{1}$ whose boundary components (diffeomorphic to circles) are also boundary 
components of $\partial \w U^{\e}$. In particular, since $\w U^{\e}$ has a bounded number of components, so does $\w U_{\e}$. 
By \eqref{dinf}, there must be some components $A^{j}$ of $\w U_{\e}$ with $diam_{g_{i}}A^{j} \to \infty$ as $i \to \infty$. On 
those components $A^{j'}$ with uniformly bounded diameter, the injectivity radius is uniformly bounded below (again by standard 
comparison geometry, cf.~\cite{Pet}). Thus, by choosing $\e$ smaller, we may assume the decomposition \eqref{tt} satisfies: $\w U^{\e}$ 
consists of a fixed number of components of uniformly bounded diameter, $\w U_{\e}$ consists of a fixed number of annuli $A^{j}$ with 
\be \label{dinf2}
diam_{g_{i}}A^{j} \to \infty  \ \ {\rm as} \ \ i \to \infty,
\ee
and so with $inj_{g_{i}}(y_{i}) \to 0$ at some points $y_{i} \in A^{j}$ between the boundary components. These thin 
regions, also frequently called cusp regions, join the thick components together to form $S^{2}$. 

  Next we claim that for any $j$, $\inf_{A^{j}}K_{i} \to 0$ as $i \to \infty$. In fact let $A_{L} \subset A^{j}$ be any subannulus 
for which the $g_{i}$-distance to its two boundary components $\partial A_{L}$ is at least $L$; then 
\be \label{neg}
\inf_{A_{L}}K_{i} \leq \pi L^{-1} .
\ee
The estimate \eqref{neg} is an immediate consequence of the standard Myers' theorem in comparison geometry (cf.~\cite{Pet}).

   The extremal length $\ell_{ext}(A)$ of an annulus $A$, cf.~\cite{AS}, \cite{Ka} for instance, is defined to be the supremum of the ratio 
$\ell^2(\s)/area(A)$, where $\ell(\s)$ is the minimal length of curves $\s$ joining distinct boundary components of $A$ and the 
supremum is over all conformally equivalent metrics on $A$. By construction, the extremal length is a conformal invariant 
of $A$ and is increasing under inclusion. By the uniformization theorem, any annulus $A$ is conformally equivalent 
to a standard annulus $A(r_1, r_2)$ of inner and outer radii $r_1$, $r_2$ in the Euclidean plane $\bR^2$. One has 
\be \label{extell}
\ell_{ext}(A(r_1, r_2)) = \frac{1}{2\pi}\log(\frac{r_2}{r_1}),
\ee
cf.~again \cite{AS}, \cite{Ka}. A similar formula holds for annuli in spherical geometry. By construction, the thin annuli $(A^{j}, g_{i})$ 
above satisfy 
$$\ell_{ext}(A^{j}, g_i) \to \infty,$$ 
as $i \to \infty$. Since $area_{g_{i}}A^j$ is bounded above by \eqref{abound}, it follows that there are curves $\g^j \subset A^j$, homologous 
to a boundary component of $A^j$, such that $\ell_{g_{i}}(\g^j) \to 0$ as $i \to \infty$. Since then $inj_{g_{i}}(q_j) \to 0$ for any point 
$q_j \in \g^j$, one has $dist_{g_{i}}(\g^{j}, \partial A^j) \to \infty$ as $i \to \infty$. For convenience, we choose $\g^j$ to be a shortest 
closed geodesic homotopic to $S^1$ in $A^j \simeq I\times S^1$. 

  Let $B^j$ be one component of $A^j \setminus \g^j$, so $B^j$ is a subannulus of $A^j$, still with $\ell_{ext}(B^j) \to \infty$. The 
thick boundary component of $B^j$ is in $\w U^\e$ while the inner or thin boundary $\g^j$ has $\ell_{g_{i}}(\g^j) \to 0$. 
In addition, we choose an annulus $A_{L} \subset B^j$ as in \eqref{neg} with $\g^j$ equal to the inner boundary component of 
$\partial A_L$.

  Up to this point, the annuli $(A^{j}, g_{i})$ and half-annuli $(B^j, g_i)$ have been considered as abstract Riemannian manifolds. 
However, since $g_{i}$ is pointwise conformal to $g_{+1}$, the annuli $(B^{j}, g_{i})$ are domains in the standard round sphere 
$S^{2}(1)$. In spherical geometry, as in \eqref{extell} the extremal length of an annulus can become unbounded only if at least 
one boundary component converges to a point, so that we may assume $\ell_{g_{+1}}(\g^j) \to 0$. In particular one may choose 
a point $p^j \in B^j$ such that $\g^j$ converges to $p^j$, in that $dist_{+1}(\g^j, p^j) \to 0$. 

  Next, the short curve $\g^j$ bounds a disc (on each side) in $S^2$ and at least one such disc $D_j^2$ is small, i.e.~contained in a spherical 
disc of (arbitrarily) small radius. Choose the component $B^j$ so that $B^j \subset D_j^2$ with $\partial D_j^2 = \g^j$. Applying the 
Gauss-Bonnet theorem to such discs we obtain 
\be \label{gb}
\int_{D_j^{2}}K_i dv_{g_i} = 2\pi - \int_{\partial D_j^{2}}\k_i,
\ee
where $\k_i$ is the geodesic curvature of the boundary $(\partial D_j^2, g_i)$. Since the boundary of $D_j^2$ is the geodesic $\g^j$, this 
gives 
$$\int_{D_j^{2}}K_i dv_{g_i} = 2\pi.$$
The area bound \eqref{abound} then implies there exists $\k_{0}$, depending only on $A_{0}$ in \eqref{abound} (and $\e$) such that 
\be \label{kap0}
K(x) \geq \k_{0},
\ee
for some $x \in D_j^{2}$. (Simple examples, e.g.~a sequence of round spheres with radius diverging to infinity, show 
that upper area bound is necessary for \eqref{kap0}). 

  We now obtain a contradiction as follows. The discs $D_j^2$ (which depend on $i$) are contained in $g_{+1}$-geodesic discs 
$D_{p^j}(\d_{i})$ with $\d_i \to 0$ as $i \to \infty$. Such discs thus contain points where both \eqref{kap0} and \eqref{neg} hold 
with $L^{-1}$ arbitrarily small. Since the sequence $\{K_{i}\}$ is uniformly bounded in $C^{\a}$, this gives a contradiction.  

{\endproof} 

  Next we turn to the existence of an apriori area bound.  Of course if $\min K_i \geq \k_0$ for some $\k_0 > 0$, then the 
area bound is automatic; the arguments below through and including Proposition 2.6 are only needed when $\liminf K_i \leq 0$. 
We first collect several general facts about the behavior of $\{u_i\}$. 
\begin{lemma}
Let $g_i = e^{2u_i}g_{+1}$ and suppose $K_i \to K \in C^{m-2,\a}$. Then the following results hold. 

{\bf (i).} For any $\e > 0$ there is a constant $A_1 = A_1(\e) < \infty$ such that 
\be \label{abound1.5}
area_{g_{i}}\{K_{i} < -\e\} \leq A_1.
\ee

{\bf (ii).} There is a constant $M = M(\e) < \infty$ such that  
\be \label{uub}
u_i(x) \leq M,
\ee
for all $x \in \{K_i \leq -\e\}$. 

{\bf (iii).} If $area_{g_i}S^2 \to \infty$, then 
\be \label{aKbound}
area_{g_{i}}\{-\e < K_{i} < \e\} \to \infty.
\ee

{\bf (iv).} If $\liminf K_i \leq 0$ and $K \in C_+^{m-2,\a}$, then there exists a constant $m_0 < \infty$ such that on $S^2$, 
\be \label{lower}
u_i \geq -m_0.  
\ee
\end{lemma} 

{\bf Proof:} The computations and estimates to follow in the proof apply to each $u_i$, but we will generally drop the 
index $i$ from the notation; thus $u = u_i$ and $K = K_i$ in the following. The estimates obtained will all be uniform, 
i.e.~independent of $i$ large, using in particular the convergence $K_i \to K \in C^{m-2,\a}$. 

  To begin, one has 
$$area_{g}\{K < -\e\} = \int_{\{K < -\e\}}e^{2u}dv_{g_{+1}}.$$
The estimate \eqref{abound1.5} is obvious in the region where $u \leq 0$, so we concentrate in the region where $u \geq 0$. 
Without loss of generality, choose $\e$ such that all points in the interval $[-\e, -\e/2]$ are regular 
values of $K$. Let $\eta = \eta(K)$ be a smooth cutoff function of $K$ with $\eta \in [0,1]$, $\eta = 1$ on $\{K \leq -\e\}$ and 
$\eta = 0$ on $\{K \geq 0\}$. Let  $U^{0} = \{u \geq 0\}$ and 
$L_{0} = \partial U^{0} = \{u = 0\}$. (In the following, one should replace $0$ by a regular value $\d$ of $u$ arbitrarily close to $0$, 
but this will make no difference in the argument below). From \eqref{u}, one has 
$$-\e \int_{U_0}\eta^4 e^{2u} \geq \int_{U_0} \eta^4 K e^{2u} \geq - \int_{U_0}\eta^4 \D u,$$
where the integration is with respect to $dv_{g_{+1}}$. 
Applying the divergence theorem to the last term gives $-\int \eta^4 \D u = -\int u \D \eta^4 + \int_{L_0}u \partial_{N}\eta^4 - \int_{L_0}\eta^4 
\partial_N u$. The first boundary integral vanishes since $u = 0$ on $L_0$ while the second is non-negative since $N$ points out of 
$U_0$. Hence 
$$-\e\int_{U_0}\eta^4 e^{2u} \geq -\int_{U_0}u \D \eta^4.$$
Since $\D \eta^4 = 4\eta^3\D \eta + 12\eta^2|d\eta|^2$, it follows that 
$$\int_{U_0}\eta^4 e^{2u} \leq C\int_{U_0}\eta^2 u,$$
where $C$ depends only on $\e$ and the $C^2$ norm of $K$ (and the fixed choice of $\eta(K)$). By the Cauchy-Schwarz inequality, one 
has $\int_{U_0}\eta^2 u \leq \sqrt{4\pi}(\int_{U_0}\eta^4 u^2)^{1/2}$. Moreover,  $e^{2u} \geq u^2$ on $U^{0}$ and $e^{2u} > > 
u$ for $u$ large, so that the integrand on the left side of the inequality above dominates the right side integrand where $u$ is large. 
It follows that $\int_{U^{0}}\eta e^{2u} \leq C'$ so that 
$$area_{g}(\{K \leq -\e\} \cap U^{0}) \leq \int_{U^{0}}\eta^4 e^{2u} \leq C'.$$
This proves \eqref{abound1.5}. 

  To prove (ii), the area bound \eqref{abound1.5} gives a uniform $L^1$ bound on $e^{2u}$, $u = u_i$, and hence a uniform $L^2$ bound 
on $u$ on $\{K \leq -\e\}$. Since on this set $\D u \geq 1$ by \eqref{u}, the upper bound \eqref{uub} is an immediate consequence of the 
well-known DeGiorgi-Nash-Moser interior estimate for subsolutions of elliptic equations, cf.~\cite{GT}. Of course here we use the estimate on 
the domains $V' = \{K \leq -\e\} \subset V = \{K \leq -\e/2\}$ for regular values $\e/2, \e$ as well as the convergence $K_i \to K \in C^{m-2,\a}$. 

  For the third claim, by the Gauss-Bonnet theorem, one has 
\be \label{gb2}
\int_{\{K \leq -\e\}}Ke^{2u} + \int_{\{K \in [-\e, \e]\}}Ke^{2u} + \int_{\{K \geq \e\}}Ke^{2u} = 4\pi,
\ee
Now \eqref{abound1.5} implies the first integral in \eqref{gb2} is bounded below. If $area_{g_{i}}\{-\e < K_{i} < \e\}$ is bounded, then 
the second integral in \eqref{gb2} is bounded above and hence so is the third integral. That in turn implies $area_{g_{i}}\{K_{i} \geq \e\}$ is 
uniformly bounded, which implies $area_{g_{i}}S^{2}$ is uniformly bounded. Hence if $area_{g_i}S^2 \to \infty$, then 
\eqref{aKbound} must hold. 

   To prove (iv), a simple computation using \eqref{u} gives 
\be \label{eu}
(\D +2)e^{-2u} = 2K + 4e^{-2u}|du|^2.
\ee
Now the first eigenvalue $\l_1$ of the Laplacian on $(S^2, g_{+1})$ is $\l_1 = 2$, with eigenfunctions given by restrictions of the linear functions 
$\ell$ on $\bR^3$ to $S^2(1)$. Any proper smooth domain $U \subset \subset S^2$ thus has lowest eigenvalue $\l_1$ for the Laplacian with 
Dirichlet boundary values on $\partial U$ satisfying $\l_1 > 2$. The operator $\D + 2$ thus has no kernel on smooth functions on $U$ 
vanishing at $\partial U$. Observe that \eqref{eu} implies $(\D + 2)e^{-2u} \geq - C$, where $C$ is a lower bound for $2K$. Let 
$U$ be any smooth domain for which there is a spherical disc $D_p(\e) \subset S^2  \setminus U$. It then follows from the well-known 
maximum principle for subsolutions of $\D + 2$, cf.~\cite{GT}, that 
\be \label{supes}
\sup_{U} e^{-2u} \leq C'(\e) + \sup_{\partial U}e^{-2u},
\ee
where $C'(\e)$ depends only on $\e$. Thus if $e^{-2u}$ (i.e.~the sequence $e^{-2u_i}$) is uniformly bounded above on $\partial U$, 
then $e^{-2u}$ is bounded, i.e.~$u$ is uniformly bounded below, on $U$, provided the complement $U^c$ contains a small disc 
$D_p(\e)$ of fixed radius. This implies that if $u_i \to -\infty$ anywhere on $S^2$, then $u_i \to -\infty$ at least on a dense set in 
$S^2$. Note that exactly this behavior occurs for the family of conformal dilations in \eqref{confact}, $u_i = \log \chi_{p,t_i}$ 
with $t_i \to \infty$. 

  Choose a point $p_i$ realizing $\max u_i$ so that $u_i(p_i) \to \infty$. Next choose a point $q_i$ with $K_i(q_i) = 0$. Without loss 
of generality, we may assume $u_i(q_i) \to -\infty$; (otherwise, one may replace $q_i$ by $q_j$ with $K_i(q_j) \to 0$ as $i, j \to \infty$ 
and $u_i(q_j) \to -\infty$ as $j \to \infty$ and take a diagonal subsequence). As in the process described in \eqref{compf}, apply a 
sequence of conformal dilations $\f_i = \f_{q_i,p_i,t_i}$ to $u_i$ with source $q_i \to q$ and sink $p_i \to p$, so that $\w u_i = \f_i^*u_i$ 
satisfies 
\be \label{renorm}
\int_{D_q(\frac{1}{2})}e^{2\w u_i} dv_{+1} = 1.
\ee
Note here that for any fixed $u$, the integral $\int_{D_q(\frac{1}{2})}e^{2\f^*u} dv_{+1}$ is a continuous function of $\f \in \mathrm{Conf}(S^2)$ 
and varies between $(a\min e^{2u}, a\max e^{2u})$, $a = area D_q(\frac{1}{2})$. Since $u_i(q_i) \to -\infty$ and $u_i(p_i) \to +\infty$, the conformal dilations $\f_i$ above exist and $t_i \to +\infty$. (The choice of the constants $\frac{1}{2}$ and $1$ here can be replaced by other fixed constants). If $\inf_{D_q(\frac{1}{10})}\w u_i \to -\infty$, then one may dilate or rescale further 
so that $\hat u_i = \f_{p_i,q_i,s_i}^*\w u_i$ is bounded below in $D_q(\frac{1}{2})$ and the area with respect to $\hat u_i$ in 
\eqref{renorm} is at most $1$. Relabeling then $\hat u_i$ to $\w u_i$ if necessary, the estimate \eqref{supes} then implies that 
$\w u_i$ is bounded below globally on $S^2$. 
 
  Now apply Remark 2.2 to the sequence $(S^2, \w g_i, q_i)$, $\w g_i = e^{2\w u_i}g_{+1}$. Since $\w u_i$ is uniformly bounded 
below,  $\w g_i \geq c g_{+1}$ for some uniform constant $c > 0$ and hence the geodesic ball $B_{q_i}^{\w g_i}(\d_0) \subset D_q(\frac{1}{2})$ for 
a uniform $\d_0 > 0$. We first note that the sequence $\w g_i$ cannot collapse at $q_i$, i.e.~\eqref{nc} holds. For if collapse did occur, 
then very short $\w g_i$-geodesic loops form near $q$. Such loops bound a spherical disc $D_i \subset D_q(\frac{1}{2})$ near $q$, 
whose $g_i$-area is uniformly bounded by \eqref{renorm}. This gives the same contradiction as at end of the proof of Proposition 2.4. 

  It then follows by Remark 2.2 that the pointed sequence $(S^2, \w g_i, q_i)$,  converges in a subsequence in the based $C^{1,\a}$ topology modulo 
diffeomorphisms to a limit $(\O, \w g, q)$, i.e.~\eqref{oconv} holds for some limit function $\w u$. The dilations $\f_i$ in \eqref{toId} are 
bounded, since otherwise one would neccessarily have $\w u_i \to -\infty$ almost everywhere, i.e.~the conformal rescaling above removes the possible divergence in the 
conformal group.  It follows from \eqref{oconv}-\eqref{toId} that $\{\w u_i\}$ subconverges to a solution $\w u$ of 
$$\D \w u = 1,$$
on a domain $\O \subset S^2$. Now the metric $g = e^{2\w u}g_{+1}$ is complete and flat and hence is either the Euclidean metric $g_{Eucl}$ on 
$\bR^2 = S^2\setminus \{p\}$ or a complete flat metric on a cylinder $\bR \times S^1 \simeq S^2 \setminus \{p \cup p'\}$. However, the 
metric $g_{Eucl}$ is not pointwise conformal to $g_{+1}$, giving a contradiction in that case. (Note that the pullback of $g_{Eucl}$ to $S^2 
\setminus \{p\}$ by stereographic projection is not pointwise conformal to $g_{+1}$, although it is of course conformally equivalent to 
$g_{+1}$). If $g$ is a complete flat metric on a cylinder, then there are annuli $A_i \subset (S^2, g_i)$ which are almost flat, i.e.~$K_i \to 0$ 
on $A_i$ and with length $\ell(A_i) \to \infty$. Since $K_i \to K \in C_+^{m-2,\a}$, this gives the same contradiction as at 
the end of the proof of Proposition 2.4, again using the area bound \eqref{renorm}. 

{\endproof}
  
\begin{proposition} 
Under the assumptions of Theorem 2.3, there is a constant $A_0 < \infty$ such that \eqref{abound} holds, i.e.
\be \label{abound2}
area_{g_{i}}S^{2} \leq A_{0}.
\ee
\end{proposition}
{\bf Proof:}  The proof is again by contradiction. Thus assume 
\be \label{areainf}
area_{g_{i}}S^{2} = \int_{S^2}e^{2u_i}dv_{+1} \to \infty. 
\ee
Consider the behavior of the (signed) measures 
$$d\mu_{i} = -\D u_idv_{+1} = (K_i e^{2u_i} - 1)dv_{+1}.$$
We first analyze the situation where these measures have uniformly bounded total variation (or mass), so that there is a 
constant $M < \infty$ such that 
$$\int_{S^2}|\D u_i|dv_{+1} \leq M.$$
This bound is equivalent to 
\be \label{abscurv}
\int_{S^2}|K_i|dv_{g_{i}} = \int_{S^2}|K_i|e^{2u_i}dv_{+1} \leq M,
\ee
so there is a uniform $L^1$ bound on the total absolute curvature of $(S^2, g_i)$. 

  Let $y_i \in S^2$ be a sequence of base points; in the following it is convenient (but not necessary) to choose $y_i$ to be 
points realizing $\max K_i$ with $y_i \to y \in S^2$. Recall that since $K_i \to K \in \cC$, there is a constant $\k_0 > 0$ such that 
$K_i(y_i) \geq \k_0$ for all $i$. Hence there is an $\e_0 > 0$ and spherical disc $D_{y_i}(\e_0)$ about $y_i$ of $g_{+1}$-radius $\e_0$ 
such that $K_i(x_i) \geq \k_0 / 2$, for $x_i \in D_{y_i}(\e_0)$. By \eqref{lower}, there is a $g_i$-geodesic disc $B_i = B_{y_i}(\d_0)$ 
of $g_i$-radius $\d_0$ about $y_i$ such that $B_{y_i}(\d_0) \subset D_{y_i}(\e_0)$. Now it is well-known that geodesic discs of 
fixed radius cannot collapse, so $inj_{g_{i}}(y_i) \geq i_0 > 0$. This follows for instance from the Gauss-Bonnet theorem:
$$\int_{B_i}K_i = 2\pi\chi(B_i) - \int_{\partial B_i}\k_i.$$
If $B_i$ collapses, then topologically $B_i$ is an annulus so $\chi(B_i) = 0$. Standard comparison geometry for the exponential map 
in local spaces with $0 \leq K \leq K_0$ shows that the geodesic curvature $\k_i > 0$, which gives a contradiction. 

  Remark 2.2 then implies that the sequence $(S^2, g_{i}, y_{i})$ converges (in a subsequence) in the pointed $C^{1,\a}$ topology 
modulo diffeomorphisms, to a complete Riemannian surface $(\O, \w g, \w y)$. The lower bound \eqref{lower}  and the fact that 
$K_i \to K$ in $C^{m-2,\a}$ implies that $|\nabla K_i|_{g_i}$ is uniformly bounded, so the convergence is actually in the 
$C^{2,\a}$ topology, again by Remark 2.2. The surface $(\O, g)$ is necessarily open; if $\O$ is closed, then $\O = S^2$ and 
$area_{g_i} S^2$ is uniformly bounded, contradicting \eqref{areainf}. 

  By \eqref{abscurv}, $(\O, \w g)$ has finite total absolute curvature. A well-known theorem of Huber \cite{Hu} then implies that ($\O, \w g)$ 
is conformally equivalent to $S^2$ punctured at a finite number of points $q_j$, i.e.~each end $E_j$ of $\O$ is parabolic in the sense 
of potential theory, cf.~\cite{AS}. We note also that for any end $E$ of $\O$, there 
exists a divergent sequence of points $x_i \to \infty$ such that 
\be \label{Kto0}
\lim_{i\to \infty} K(x_i) \leq 0.
\ee
(For if \eqref{Kto0} were not the case, then $K \geq \k > 0$ on $E$, for some constant $\k > 0$, which implies that $E$ is of bounded 
diameter by the Myers theorem in comparison geometry, \cite{Pe}; this is of course a contradiction). 

  Referring again to Remark 2.2, without loss of generality, we may assume $\psi_i^{-1}(y_i) = \w y_i \to \w y$ with $F(\w y) = y$. 
Since $K_i(y_i) > 0$ is bounded away from 0, it follows from \eqref{Kto0} and the $C^{\a}$ convergence $K_i \to K$, that the 
image of the embeddings $\psi_i$ must contain a small spherical disc of fixed size, i.e.~$D_{y_i}(\d) \subset \psi_i(\O_i)$ 
for some $\d > 0$. The embedding $F: \O \to S^2$ has image $F(\O) = S^2$ with a finite number of points $\{q_j\}$ removed. 
As in \eqref{oconv} one has 
\be \label{fsu}
(F^{-1})^*\psi_i^*(e^{2u_i}g_{+1}) \to e^{2u}g_{+1},
\ee
on $F(\O)$ and there exists a sequence of conformal transformations $\f_i$ such that 
$$\f_i \circ \psi_i \circ F^{-1} \to Id,$$
on $S^2$. Now we claim that the sequence $\f_i$ cannot diverge to $\infty$ in $\mathrm{Conf}(S^2)$. Namely, since the images 
$\psi_i(\O_i)$ contain a spherical disc of fixed size about $y$, $\f_i$ can only diverge if it is within a bounded distance in 
$\mathrm{Conf}(S^2)$ to a sequence of conformal dilations $\bar \f_i = \f_{y,t_i}$ sourced at $y$ with $t_i \to \infty$. Setting 
$\w u_i = \bar \f_i^*u_i$, it follows from \eqref{fsu} that $\w u_i$ subconverges in $C^{\a}\cap W_{*}^{1,2}$ to a solution 
$\w u$ of $\D \w u = 1 - K(y)e^{2\w u}$. By elliptic regularity, $\w u$ is smooth. This is the equation for a metric of constant 
curvature $K(y)$ which implies that $\O$ is compact, i.e.~$\O = S^2$, giving a contradiction since $\O$ is open. 

 Thus $F^{-1}\circ \psi_i$ converges to a limit quasi-conformal diffeomorphism with dilatation $1$, i.e.~to a conformal 
diffeomorphism. Altering the embedding $F$ by this fixed conformal map, we may assume $F$ is the inclusion of a subset $\O \subset S^2$ 
and 
$$\O = S^2 \setminus \cup \{q_j\}.$$
By \eqref{oconv}, the functions $u_i$ converge to the limit function $u$ on $\O$. Elliptic regularity associated with the equation 
\eqref{u} implies the convergence is in $C^{m,\a}$. Since $K = \lim K_i$, by \eqref{Kto0} one must have $K(q_j) \leq 0$. 
In fact, one must have 
$$K(q_j) = 0.$$ 
Namely if $K(q_j) < 0$, then the end $E_j$ of $\O$ is a hyperbolic cusp and so of finite area. 
The end $E_j$ must be capped off by a small spherical disc $D_j$ about $q_j$ in $S^2$. Applying the Gauss-Bonnet theorem to 
$(D_j, g_i)$ for $i$ large shows that there exist points $q_j' \in D_j$ such that $K_i(q_j') \geq 0$. For $i$ sufficiently large, this 
gives the same contradiction as before, namely to the $C^{\a}$ convergence $K_i \to K$. 

   Now the measures $d\mu_i$ converge (in a subsequence) to a limit signed measure $d\mu$ on $S^2$. Let $d\mu_r$ and 
$d\mu_s$ denote the absolutely continuous and singular measures of $d\mu$ with respect to Lebesque measure:
\be \label{regsing}
d\mu = d\mu_r + d\mu_s.
\ee
As above, the functions $u_i$ converge (in a subsequence) uniformly on compact 
subsets of $\O$ in $C^{m,\a}$, to a limit function $u$. The regular part of $d\mu$ is thus given by 
$$d\mu_r = (Ke^{2u} - 1)dv_{+1} \ \ {\rm on} \ \ \O,$$
where $K = \lim K_i$. Since $area_{g_{i}}S^2 \to \infty$, $u_i$ blows up, $u_i \to \infty$, at each 
$q_j$. Similarly, since $g$ is complete, $u \to \infty$ at each $q_j$. As noted above, the support of the singular measure $d\mu_s$ is given by 
$$Z = \cup q_j \subset \{K = 0\}.$$
Any measure supported at a point is a multiple of the Dirac delta measure, so that there exist $a_j \in \bR$ such that 
\be \label{sing}
d\mu_s = 2\pi \sum a_j \d_{q_{j}}.
\ee
One has 
$$\int_{\O}d\mu_r = \int_{\O}Kdv_g - 4\pi \leq 2\pi \chi(\O) - 4\pi,$$  
where the last inequality follows from another well-known result of Huber \cite{Hu}. Since the total mass of $d\mu$ is $0$ 
(by the weak convergence of the measures $d\mu_i \to d\mu$), it follows that 
$$\int d\mu_s = 2\pi \sum a_j \geq 4\pi - 2\pi \chi(\O),$$
so that $\sum a_j \geq 2 - \chi(\O)$. Note that the Euler characteristic $\chi(\O)$ is given by $\chi(\O) = 2 - \#\{\cup q_j\}$.

Let $G = G_q$ be the Green's function for the operator $1 - \D$ on $S^2(1)$ with pole at $q$ (and $-q$), so that 
$$\D G =  1 - 2\pi (\d_q + \d_{-q}).$$
Thus $G_q$ corresponds to the Dirac measure $\d_{q}$ with $a = 1$ (and measure $\d_{-q}$). The behavior of the 
Green's function at the antipodal point $-q$ will play no role in the discussion below. A simple computation shows that 
$G$ is given explicitly by 
\be \label{Greenlog}
G(x,q) = \log \frac{1}{\sin r(x,q)}
\ee
where $r(x,q) = dist_{+1}(x,q)$. In particular for $x$ near $q$, $G_q$ has the well-known logarithmic behavior $G(x, q) = -\log r(x, q) + 
\a$, where $\a$ is smooth. 

  Now we first note that \eqref{Greenlog} implies that in \eqref{sing}
$$a_j \geq 1.$$
This follows from the completeness of $(\O, g)$. Namely, the $g$-length of a curve $\g$ ending at some $q_j$ is given by 
\be \label{length}
\ell_g (\g) = \int_{\g} e^u dr,
\ee
where $r$ is the $g_{+1}$-arclength parameter for $\g$. Since $g$ is complete, the integral in \eqref{length} must diverge; hence 
on approach to $q_j$, $u \geq \log r^{-1+\d}$, for any fixed $\d > 0$ (on some sequence of points $p_k$ with $r(p_k) \to 0$). 
By \eqref{Greenlog}, this implies $a_j \geq 1$. 

  The leading order behavior of the limit function $u$ near the point $q \in \{q_j\}$ is governed by the behavior of the Green's function 
$G = G_q$. By Green's representation formula, the decomposition \eqref{regsing} gives the expression 
\be \label{uGb}
u = aG + \b,
\ee
near $q$ where $\b$ is (at least) continuous in a neighborhood of $q \in \{q_j\}$.  This gives 
$$g = e^{2u}g_{+1} = e^{2(aG + \b)}g_{+1} = e^{2\b}g_a,$$ 
where $g_a = e^{2aG}g_{+1}$. One has 
$$K_a e^{2aG} = 1 - \D aG = (1-a) + 2\pi a(\d_q + \d_{-q}),$$
so that away from the pole $q$ (and $-q$), 
\be \label{ka}
K_a = (1-a)e^{-2aG} = (1-a)(\sin r)^{2a} \leq 0 .
\ee
For $K_u$ we obtain 
\be \label{Ku}
K_u = (1-a - \D \b)(\sin r)^{2a} e^{2\b} = \chi r^{2a},
\ee
on $\O$ near $q$; the last equality defines $\chi$. Now $K_u = \lim K_i = K$, so that \eqref{Ku} holds also at $q$. 
Clearly then $K(q) = 0$ and $q$ is a critical point of $K \in \cC$. Since \eqref{criterion} thus holds at $q$, one must have 
$\D K(q) > 0$. 
This implies that 
\be \label{a=1}
a = 1,
\ee
and also $\chi(q) > 0$. 

   Now the geometry of $g_a$ for $a = 1$ is that of a complete flat metric on the cylinder $\bR\times S^1$; $K_a = 0$ when $a = 1$. 
The pole $q$ occurs at $+\infty$ while the pole $-q$ occurs at $-\infty$. As above, we ignore the behavior at $-q$. It follows 
that the geometry of $(\O, g)$ on approach to any $q \in \{q_j\}$ is that of a flat cylindrical end. 

   On the sequence $(S^2, g_i)$, these cylindrical ends are filled in or capped by very small discs $D_j = D_j(\d_i)$ in $(S^2, g_{+1})$ 
about $q_j$, $\d_i \to 0$. Thus topologically, $S^2 = \O \cup \{\cup D_j\}$. The metrics $g_i$ are blowing up such very small discs either 
to discs of approximately unit area or to very large discs. The former case corresponds to capping of the cylindrical end $\bR^{+}\times 
S^1$ by a hemisphere. In this case, the Gauss curvature is necessarily bounded away from $0$ at some point in the disc, (by the 
Gauss-Bonnet theorem for instance), which contradicts the $C^{\a}$ convergence $K_i \to K$. Hence the capping discs must become 
very large and of very small curvature $K_i \sim 0$ corresponding to a spherical bubble with area tending to infinity. Intuitively, this 
requires a band of negative curvature to open up the cylinder, capped off by a region of positive curvature to cap off the end. 

  To analyse this latter case in detail, fix any $i$ sufficiently large; for $r$ very small, by \eqref{Ku} the curvature $K_i$ of $g_i$ on the sphere 
$S_q(r)$ is small but uniformly positive, $K_i \geq \k_0 > 0$. Consider the exponential map of $g_i$-geodesics normal to $S_q(r)$ into the 
disc $D_q(r)$. If the curvature $K_i$ of $g_i$ satisfies $K_i \geq 0$ in $D_q(r)$, then standard comparison geometry implies these geodesics 
have conjugate points at a $g_i$-bounded distance to $S_q(r)$. This implies $(D_q(r), g_i)$ has uniformly bounded diameter, which 
contradicts the assumption here that $area_{g_i}D_q(r) \to \infty$. Thus, $\min_{D_q(r)} K_i < 0$ on some sequence $i \to \infty$. 
Let $L_0 = L_0(i) \subset \{K_i = 0\} \subset D_q(r)$, be the collection of curves homologous to $S_q(r)$ such that $K_i \geq 0$ in the 
region bounded by $S_q(r) \cup L_0$. Now on the one hand, standard comparison geometry for non-negative curvature implies the 
$g_i$-length $\ell_{g_i}$ of $L_0$ is uniformly bounded. Each component $L_0'$ of $L_0$ bounds a small disc $D_i' \subset D_q(r)$ in $S^2$. 
On the other hand, standard comparison geometry for non-positive curvature shows that if $K_i \leq 0$ on $D_i'$, then $area_{g_i} D_i' \to \infty$ 
implies $\ell_{g_i}(L_0') \to \infty$; this is a  standard isoperimetric-type inequality for non-positive curvature. Since $area_{g_i}D_q(r) \to 
\infty$, it follows that for some $D_i'$, $\max_{D_i'} K_i > 0$ for $i$ large. However, this contradicts the assumption that $K_i \to K \in \cC$. 
Namely, since $\{K_i\}$ is uniformly within $\cC$, there are no local maxima of $K_i$ with $|K_i|$ sufficiently small. 

It is worth noting, and is discussed in detail in Remark 2.7 below, that when $K \notin \cC$, such bubble behavior with unbounded 
area may occur. 

\medskip 

   Next suppose the total mass or variation of $d\mu_i$ satisfies
$$m(|d\mu_i|) \to \infty.$$
In this case, we will see that the blow-up behavior of $\{u_i\}$ creates even larger $a$ in \eqref{Ku}, leading to the same contradiction.  

   To begin, as above choose base points $y_i$ realizing $\max K_i$. By Proposition 2.1 and Remark 2.2 as above, the 
sequence $(S^2, g_i, y_i)$ subconverges to a complete conformally flat limit surface $(\O, g, y)$; the convergence is uniform in 
$C^{2,\a}$ (at least) on compact subsets of $\O$. As in Remark 2.2, the abstract domain $(\O, g)$ may be conformally embedded 
as a planar domain, i.e.~a domain in $S^2$ with complete metric $g = e^{2u}g_{+1}$. As before we let $Z = S^2 \setminus \O$ be the 
singular set of $g$. By \eqref{areainf}, there is a sequence of nested neighborhoods $N_m \supset N_{m+1}$ of $Z$ in $S^2$ with 
$\cap N_m = Z$ such that $area_{g_i}N_m \to \infty$ for each fixed $m$. By completeness $u \to \infty$ at $\partial \O = \partial Z$, 
and hence by \eqref{uub},
\be \label{Z0}
\partial Z \subset \{K \geq 0\},
\ee
where $K = \lim K_i$. (It will follow from later arguments below that $Z$ has empty interior, so that $\partial Z = Z$). 

   The complete surface $(\O, g)$ may have a finite or infinite number of ends. Suppose first the number of ends is finite and 
choose any fixed end $E$. Then $E$ is an annulus in $S^2$, with either finite or infinite extremal length. 
This dichotomy corresponds to a hyperbolic or parabolic end respectively, in the sense of potential theory, cf.~\cite{AS}. 

  We claim that the end must be parabolic. This is proved by contradiction, so suppose $E$ is hyperbolic. Let $\s_j$ be 
a divergent sequence of smooth closed curves in the annular end $E$ disconnecting $E$. The annulus $E$ is embedded in 
$S^2$ and so $\{\s_j\}$ is a sequence of Jordan curves in $S^2$ bounding a nested sequence of domains $D_{j+1} 
\subset D_j \subset S^2$ with $\s_k \subset D_j$ for $k > j$. The limit $D_{\infty} = \lim D_j = \cap D_j$ is a 
subset of $Z$ in \eqref{Z0}. Since the end is hyperbolic, $D_{\infty}$ cannot be a point (since $E$ is parabolic 
if $D_{\infty}$ is a point), and hence $\partial D_{\infty}$ contains a non-trivial arc, i.e.~the $C^0$ image of an interval. 

    Let $r$ be the $g$-geodesic distance from a fixed point $p_0 \in E$. By a theorem of Doyle-Grigorian, cf.~\cite{D}, \cite{G}, 
there is a divergent sequence of points $p_j \in E$, $r(p_j) \to \infty$, such that $K(p_j) \leq - (1/r^2\log r)(p_j)$. Passing 
to a subsequence, we may assume $p_j \to q$, for some $q \in \partial Z$. If $K(q) < 0$, then one has a contradiction to \eqref{Z0}, 
so that $K(q) = 0$. If $K \leq 0$ in a neighborhood of $q$, then $0$ is a local maximum of $K$, contradicting $K \in \cC$. Hence there 
exist points arbitrarily near $q$ where $K > 0$ and so the level set $K = 0$ containing $q$ locally disconnects $S^2$ into 
regions where $K > 0$ and $K < 0$. The closed curves $\s_j$ then contain arcs $\a_j^{\pm} \subset \s_j$ with $\a_j^+ \subset \{K > 0\}$ 
while $\a_j^- \subset \{K < 0\}$ with $\a_j^{\pm}$ converging to arcs $\a^{\pm} \subset \{K = 0\}$. However, there are then 
long triangles $\D_j$ with vertices $\{p_0$, $\partial(\a_j^+)\}$ with $K > 0$ in the region in $\D_j$ where $r$ is sufficiently large. 
Such triangles have extremal length diverging to $\infty$ (positive curvature gives parabolic ends) and hence $\a^+$ is a single point 
$q$, giving again a contradiction.  

  Thus the end $E$ is parabolic, and hence conformally equivalent to a punctured disc $D^2\setminus q$. As above, we work 
with the positive measures $|d\mu_i| = |\D u_i|dv_{+1}$ restricted to the end $E$, i.e.~a punctured neighborhood of $q$. 
Without loss of generality, assume $m_i = m(|d\mu_i|) \to \infty$. Consider the rescaled measures $d\w \mu_i = 
m_i^{-1}d\mu_i$, so the total mass of $d \w \mu_i$ is $1$. It follows that the measures $|d\w \mu_i|$ subconverge to a 
positive measure supported at the point $q$, i.e.~to the Dirac measure $\d_q$. Hence,  to leading order, 
$d\mu_i = m_i \d_{q} +d\nu_i$ where $n_i = m(|d\nu_i|) < < m_i$. If $n_i \to \infty$, one may repeat this process with 
$d\mu_i - m_i \d_q$ in place of $d\mu_i$ and conclude that 
$$d\mu_i = (1 - \D u_i)dv_{+1} = m_i' \d_{q} +d\nu_i$$
with $n_i = m(|d\nu_i|) \leq N$, for some $N < \infty$. As in \eqref{uGb}, this implies that on $E$, 
\be \label{uexp}
u_i = m_i' \log r^{-1} + \b_i
\ee
where $\b_i$ is bounded near $q$. Since $K_i$ is thus of the form \eqref{Ku} with $a = m_i' > 1$, this contradicts the 
fact that $K_i \to K \in \cC$. 

  We note that the same argument applies if there is at least one end of $\O$ of finite topological type. 
  
    Suppose finally all ends are of infinite topological type; in the following we work with just one end $E$. Such planar ends 
$E$ have the topological form of an infinite union of pairs of pants $P_k$; $P_k$ is diffeomorphic to $S^2$ with three discs 
removed.  Each leg of $P_k$ is glued to the waist of one of the two successors $P_{k+1}$ or $P_{k+1}'$ of $P_k$. Such ends 
have singular set $Z$ equal to a Cantor set in $S^2$, cf.~\cite{Ri} for a detailed discussion. 
 
  Choose a base point $p_0 \in E$ and let $\s$ be any geodesic ray from $p_0$ diverging to infinity in $E$. Let $\{P_k'\} \subset E$ 
be the collection of pairs of pants intersecting $\s$. This forms a (connected) sub-end $E' \subset E$ with limit given by a single 
point $q \in Z$. Next we cap off all boundary components, i.e.~circles, of $\{P_k'\}$ not 
intersecting $\s$ by discs $D_{k'}$, giving then an annular end $\w E \simeq \s \times \bR^+$ ending at $q$. 
For any fixed disc $D_{k'}$, the functions $u_i$ are unbounded on $D_{k'}$, (since $dist_{g_{i}}(D_{k'}, Z) \to \infty$), but 
are uniformly bounded on the boundaries $\partial D_{k'} \simeq S^1$. For $i$ sufficiently large, (depending on 
$k'$), we then replace $u_i|_{D_{k'}}$ by the harmonic extension $\w u_i$ of the boundary values of $u_i$ on $\partial D_{k'}$, 
i.e.~the solution of the Dirichlet problem $\D \w u_i = 0$ on $D_{k'}$ with boundary data $u_i$ on $\partial D_{k'}$. In 
particular, the maximum and minimum values of $\w u_i$ on $D_{k'}$ are bounded by the maximum and minimum values of $\w u_i$ on 
$\partial D_{k'}$. 

  This construction gives a parabolic annular end $(\w E, \w g_i)$, $\w g_i = e^{2\w u_i}g_{+1}$ limiting at a singular point 
$q \in Z$. Suppose first the mass 
$\w m_i$ of $|d\w \mu_i|$ remains uniformly bounded. Then the same analysis as that leading to \eqref{Ku} gives  
$$\w K = (1-a-\D \b)(\sin r)^{2a}e^{2\b}.$$
Of course $\w K \neq K$, but $\w K = K$ in a neighborhood of the geodesic ray $\s$, i.e.~in the neighborhood $\cup P_k'$ 
of $q$. Since $\s$ may be chosen to be any geodesic ray and since $K \in C^{m-2,\a}$, it follows that $K(q) = 0$ and $K$ 
is of the form \eqref{criterion} near $q$. This brings one exactly to the situation above where \eqref{Ku}-\eqref{a=1} holds and 
the same arguments as before lead to the same contradiction. Finally, the same argument as above concerning \eqref{uexp} also 
gives the contradiction when $\w m_i \to \infty$. This completes the proof of Proposition 2.6.  
  
{\endproof}  

  Propositions 2.4 and 2.6 together complete the proof of Theorem 2.3. 

\begin{remark}
{\rm It is important to observe that the condition \eqref{maxK} that $K \in \cC$ is necessary in Theorem 2.3. This follows 
from the construction of large spherical ``bubbles" in a fixed conformal class by Ding-Liu in \cite{DL} and later by 
Borer-Galimberti-Struwe in \cite{BGS}. In particular, the construction in \cite{BGS} shows that there is a sequence 
$\d_{i} \to 0$ and conformal metrics $g_{i} = e^{2u_{i}}g_{+1}$ on small spherical discs $D_{p}(\d_{i}) \subset S^{2}(1)$ 
with arbitrarily large area, with arbitrarily small but positive Gauss curvature $K_{i} > 0$, and for which $u_{i}$ is smoothly bounded 
near the boundary of the fixed disc $\partial D_{p}(\d)$. In particular $u_{i} \to \infty$ in $D_{p}(\d_{i})$. The construction 
has the property that the Gauss curvature $K_{i}$ has a local maximum at $p$ with $0 < K_{i}(p) = \e_{i} \to 0$ but 
$K_{i} < 0$ and small near $\partial D_{p}(\d)$ for $\d > 0$ small. Further the functions $u_{i}$ are smoothly 
bounded near $\partial D_{p}(\d)$. The functions $u_{i}$ and the associated metrics $g_{i}$  may then be extended 
past $\partial D_{p}(\d)$ to smoothly bounded conformal metrics on all of $S^{2}$. 

   This construction shows that $K_{i} \to K$ in $C_{+}^{m-2,\a}$ on $S^{2}$ does not imply that conformal metrics 
$g_{i}$ are uniformly bounded modulo diffeomorphisms, i.e.~the condition $K \in \cC$ is necessary. Clearly such bubble 
formation could occur at any assigned finite number of points in $S^{2}$ for instance. 

}
\end{remark}

  Theorem 2.3 implies that the hypotheses \eqref{ad} of Proposition 2.1 hold, for $g_{i} = e^{2u_{i}}g_{+1} \in [g_{+1}]$ 
with $K_{i} \to K \in \cC$. Proposition 2.1 then shows that convergence of the curvatures implies convergence 
of the metrics $g_{i}$ in $C^{1,\a}\cap W_{*}^{2,p}$, modulo diffeomorphism. Next we turn to the role of the 
diffeomorphisms. 

\medskip 

  As discussed following the proof of Proposition 2.1, there are diffeomorphisms $\psi_i$ of $S^2$ such 
that,  (in a subsequence) 
\be \label{a3}
\psi_{i}^{*}(e^{2u_{i}}g_{+1}) \to e^{2u}g_{+1},
\ee
in $C^{1,\a} \cap W_{*}^{2,p}$. Further, as in \eqref{id1}, there are $\f_i \in \mathrm{Conf}(S^2)$ such that 
\be \label{id}
\eta_{i} : = \f_{i}^{-1}\circ \psi_{i} \to Id.
\ee
The diffeomorphisms $\psi_i$, $\eta_{i}$ (and their inverses) are $(1+\e_{i})$ quasi-conformal 
with $\e_{i} \to 0$ as $i \to \infty$. 

The next main step is the following: 

\begin{proposition} 
For $g_{i} = e^{2u_{i}}g_{+1}$ as above with $K_{i} \to K \in \cC \subset C_{+}^{m-2,\a}$, the functions $\f_{i}^{*}u_{i}$ 
are uniformly bounded in $C^{1,\a} \cap W^{2,p}$ for any $\a < 1$, $p < \infty$. Thus $\{u_{i}\}$ is bounded in 
$C^{1,\a}\cap W^{2,p}$ modulo the action of the conformal group on $C^{m,\a}$. 
\end{proposition}

{\bf Proof:} By the uniform H\"older regularity of locally bounded quasi-conformal maps, 
cf.~\cite{LV}, \cite{GT}, the convergence in \eqref{id} is in $W_{*}^{1,2}$ and $C^{\a}$ for any fixed $\a < 1$. 
More significantly, an important theorem of Astala \cite{As} on the $W^{1,p}$ boundedness of $K$-quasi-conformal 
mappings shows that $\eta_{i}$ in \eqref{id} are bounded in $W^{1,p}$ for any $p < \infty$, so that 
$$\eta_{i} \to Id,$$
in $W_{*}^{1,p}$, i.e.~weakly in $W^{1,p}$. The same statements hold of course for $\eta_{i}^{-1}$. 

 Write then $\psi_{i}^{*} = \eta_{i}^{*}\f_{i}^{*}$ and insert in \eqref{a3} to obtain
\be \label{a4}
\eta_{i}^{*}(e^{2\f_{i}^{*}u_{i}}g_{+1}) \to e^{2u}g_{+1},
\ee
in $C^{1,\a}$, bounded in $W^{2,p}$. 
By Astala's theorem \cite{As}, $\eta_{i}^{*}$ and $(\eta_{i}^{-1})^{*}$ are uniformly bounded in $L^{p}$, for any 
$p < \infty$, and hence $e^{2\phi_{i}^{*}u_{i}}$ is uniformly bounded in $L^{p}$. 

   Now from the defining equation \eqref{u}, we have 
\be \label{uf}
\D\f_{i}^{*}u_{i} = 1 - (K_{i}\circ \f_{i})e^{2\f_{i}^{*}u_{i}}.
\ee
The term $K_{i} \circ \f_{i}$ is bounded in $L^{\infty}$ while, as noted above, the term $e^{2\f_{i}^{*}u_{i}}$ is bounded in $L^{p}$. 
Standard elliptic regularity applied to \eqref{uf} shows that $\f_{i}^{*}u_{i}$ bounded in $C^{1,\a} \cap W^{2,p}$ modulo constants, 
i.e.~the mean value of $\f_{i}^{*}u_{i}$. Since $e^{2\f_{i}^{*}u_{i}}$ is bounded in $L^{p}$, so is the mean value of $\f_{i}^{*}u_{i}$. 
This proves the result. 
  
{\endproof}

  Given these results, we now restate and prove Theorem 1.1. As in the Introduction, let 
\be \label{cN22}
\cN = \{K \in C_{+}^{m-2,\a}: |\nabla K|(p) + |\D K|(p) > 0, \forall p \ {\rm s.t.} \ K(p) > 0 \}. 
\ee
It is easy to see that $\cN$ is open and dense in $C_{+}^{m-2,\a}$ and is invariant under the action of $\Co$. Let  
$$\cU = \pi^{-1}(\cN \cap \cC),$$
and let 
\be \label{pi01}
\pi_{0} = \pi|_{\cU}: \cU \to \cN \cap \cC,
\ee
be the restriction of $\pi$ to $\cU$. 
  
\begin{theorem}
On the domain $\cU$, the curvature map 
$$\pi_{0}: \cU \to \cN \cap \cC$$
in \eqref{pi01} is a proper Fredholm map of index 0. 
\end{theorem} 

{\bf Proof:} As noted in the Introduction, $\pi$ is a Fredholm map of index $0$ so the issue is to prove $\pi_{0}$ is proper. 
Suppose $u_{i}$ is a sequence in $\cU$ with $K_{i} = K_{e^{2u_{i}}g_{+1}} = \pi_{0}(u_{i})$ satisfying 
\be \label{kconv}
K_{i} \to K \ \ {\rm in} \ \ \cN \cap \cC \subset C_{+}^{m-2,\a}. 
\ee
Then by Proposition 2.8, there are $\f_{i} \in {\rm Conf}(S^{2})$ such that the sequence $v_{i} := \f_{i}^{*}u_{i}$ converges (in a 
subsequence) to a limit in $C^{1,\a}$ and weakly in $W^{2,p}$ for any $p < \infty$. Let $v_{0} = \lim v_{i} \in C^{1,\a} \cap W^{2,p}$. 
These functions $v_{i}$ satisfy 
\be \label{vi}
\D v_{i} = 1 - (K_{i}\circ \f_{i})e^{2v_{i}}.
\ee
By means of a further (uniformly bounded) sequence of conformal transformations if necessary, we may assume without loss of 
generality, that $v_{i} \in \cS$ for all $i$, and $v_{0} \in \cS$, for $\cS$ as in \eqref{cS}. 

  If the sequence of conformal transformations $\f_{i}$ are bounded in $\Co$, then they converge (smoothly) in a subsequence 
to a limit and hence $\{u_{i}\}$ converges in $C^{1,\a}\cap W_{*}^{2,p}$ to a limit $u$ satisfying \eqref{u} weakly. Standard 
elliptic regularity theory then implies $u \in C^{m,\a}$ and, by \eqref{kconv}, the convergence $u_{i} \to u$ is $C^{m,\a}$. 
This proves the result in this case. 

   Assume then $\f_{i} \to \infty$ in $\Co$. This occurs only if the corresponding conformal dilations $[\f_{i}] \to \infty$ in 
$\Co/O(3) = \bR^{3}$, cf.~\eqref{h3}. Without loss of generality, we may then assume that, in a subsequence, $\f_{i} = \f_{q}(t_{i})$ 
is a sequence of conformal dilations with pole at $q$ and $t_{i} \to \infty$. 
Taking the limit of \eqref{vi}, it follows that 
\be \label{v0}
\D v_{0} = 1 - K(q)e^{2v_{0}}.
\ee
Namely $e^{2v_{i}} \to e^{2v_{0}}$ in $C^{1,\a}$ while $K\circ \f_{i} \to K(q)$ strongly in $L^{p}$. It follows that $v_{0}$ 
in $W^{2,p}$ is a weak solution of \eqref{v0}. As above, elliptic regularity theory implies $v_{0}$ is smooth. The 
equation \eqref{v0} means that the metric $g_{0} = e^{2v_{0}}g_{+1}$ is of constant curvature $K \equiv K(q)$ on $S^{2}$. 
As is well-known, this means that $v_{0} = \f^{*}(const)$ and since $v_{0} \in \cS$, it follows (cf.\eqref{cslice}) that 
$v_{0} = const$ and 
\be \label{vnorm}
K(q)e^{2v_{0}} = 1.
\ee  
In particular, note that $\f_{i} \to \infty$ in $\Co$ is only possible at points $q$ where $K(q) > 0$. 

  Returning to \eqref{vi}, write \eqref{vi} as 
\be \label{Vi}
\D v_{i} = 1 - (K_{i}\circ \f_{i})e^{2v_{i}} = 1 - K(q)e^{2v_{i}} + (K(q) - K_{i}\circ \f_{i})e^{2v_{i}}.
\ee
In the following, it will simplify the argument to rescale the sequence $K_{i}$ slightly, $K_{i} \to \frac{K(q)}{K_{i}(q)}K_{i}$ 
so that 
\be \label{ki}
K_{i}(q) = K(q),
\ee
for $i$ sufficiently large. This is equivalent to shifting $v_{i}$ by small constants $v_{i} \to v_{i} + \d_{i}$ with 
$\d_{i} \to 0$ as $i \to \infty$. We assume this has been done, without changing the notation. 

Using \eqref{vnorm}, \eqref{Vi} may be rewritten as 
\be \label{vi2}
\D(v_{i} - v_{0}) = 1 - e^{2(v_{i} - v_{0})} + (E_{i} + F_{i})e^{2v_{i}},
\ee
where 
\be \label{EF}
E_{i} = (K(q) - K\circ \f_{i})e^{2v_{i}} \ \ {\rm and} \ \ F_{i} = [(K - K_{i})\circ \f_{i}]e^{2v_{i}}.
\ee   
By a careful but straightforward analysis of the asymptotic behavior of $K \circ \f_{t}$ as $t \to \infty$, 
Chang-Yang show in \cite{CY3} that 
\be \label{kest1}
|K(q) - K\circ \f_{t})|_{L^{2}}^{2} \leq \chi_{q}(t) := \left\{ \begin{array}{ll} 
                              c\frac{\log t}{t^{2}} & \mbox{ if $\nabla K(q) \neq 0$}, \\
                              \frac{c}{t^{2}} & \mbox{if $\nabla K(q) = 0$.}
                              \end{array} 
                         \right.   
\ee
Using \eqref{ki}, and replacing $K$ above by $K - K_{i}$, the same argument shows that 
\be \label{kest2}
|(K - K_{i})(q) - (K - K_{i})\circ \f_{i} |_{L^{2}}^{2} \leq \chi_{q}(t_{i}). 
\ee
Since $v_{i}$ is uniformly bounded in $C^{1,\a}$, it follows that 
\be \label{eest}
|E_{i}|_{L^{2}}^{2} + |F_{i}|_{L^{2}}^{2} \leq \chi_{q}(t_{i}).
\ee

  Next for the second term on the right in \eqref{vi2}, since $v_{i} \to v_{0}$ in $C^{1,\a}$, $e^{2(v_{i}-v_{0})} = 1 + 
2(v_{i} - v_{0}) + \s_{i}$, where 
$\s_{i} = O((v_{i} - v_{0})^{2})$. Then \eqref{vi2} gives 
$$\D(v_{i} - v_{0}) + 2(v_{i} - v_{0}) = E_{i} + F_{i} - \s_{i}.$$
The kernel $E_{2}$ of the operator $\D + 2$ consists of $1^{\rm st}$ eigenfunctions of the Laplacian on $S^{2}$, i.e.~the 
restrictions of linear functions $\ell$ to $S^{2}(1)$. The tangent space $T_{0}\cS$ to $\cS$ at $v_{0} = const$ is the orthogonal 
complement of $E_{2}$. Since $v_{i}, v_{0} \in \cS$, and since $\s_{i}$ is lower order, it follows that $v_{i} - v_{0}$ satisfies 
the elliptic estimate 
\be \label{vest}
|v_{i} - v_{0}|_{L^{2,2}}^{2} \leq \chi_{q}(t_{i}).
\ee

  On other hand, pairing \eqref{vi} with $xe^{-2(v_{i}-v_{0})}$ where $x$ is any linear function on $S^{2}(1)$ gives 
$$xe^{-2(v_{i}-v_{0})}\D v_{i} = xe^{-2(v_{i}-v_{0})} - x\frac{K_{i}\circ \f_{i}}{K(q)}.$$
Since $e^{-2(v_{i}-v_{0})}\D (v_{i} - v_{0}) = -\frac{1}{2}\D e^{-2(v_{i}-v_{0})} + 2|d v_{i}|^{2}e^{-2(v_{i}-v_{0})}$, 
we obtain  
$$-{\tfrac{1}{2}}\int x\D e^{-2(v_{i}-v_{0})} + 2\int x|dv_{i}|^{2}e^{-2(v_{i}-v_{0})} = \int xe^{-2(v_{i}-v_{0})} - 
\frac{1}{K(q)}\int x K_{i}\circ \f_{i}.$$
(Here and below, all integrals are with respect to the standard volume form of $g_{+1}$). 
Using the self-adjointness of the Laplacian $\D$ and the fact that $x$ is a first eigenfunction, one sees that the 
first and third terms above cancel, leaving 
$$2\int x|dv_{i}|^{2}e^{-2(v_{i}-v_{0})} =  - \frac{1}{K(q)}\int x K_{i}\circ \f_{i}.$$
Since $|x| \leq 1$, $v_{0} = const$ and $v_{i} - v_{0} \to 0$ in $C^{0}$ it follows from \eqref{vest} that 
$$|\int x|dv_{i}|^{2}e^{-2(v_{i}-v_{0})}| < < \chi_{q}(t_{i}),$$
and so, for any $x$,  
\be \label{small}
|\int xK_{i}\circ \f_{i}| < < \chi_{q}(t_{i}).
\ee

  Now the analysis by Chang-Yang in \cite{CY3} again shows that for $K \in \cN$, i.e.~$K$ non-degenerate in the 
sense of \cite{CY3}, and for a suitable choice of $x$, one has 
\be \label{cy2}
 |\int xK\circ \f_{i}| \geq  \psi_{q}(t_{i}),
\ee
where
\be \label{cy21}
\psi_{q}(t)    = \left\{ \begin{array}{ll} 
                              \frac{c}{t} & \mbox{ if $\nabla K(q) \neq 0$,} \\
                              c\frac{\log t}{t^{2}} & \mbox{if $\nabla K(q) = 0$, $\D K(q) \neq 0$.}
                              \end{array} \right. 
\ee
As above, writing $K_{i}\circ \f_{i} = (K_{i} - K)\circ \f_{i} + K\circ \f_{i}$, the same analysis by Chang-Yang in \cite{CY3} as above 
shows that, for any $x$,  
$$|\int x(K - K_{i})\circ \f_{i}| < <   \psi_{q}(t_{i}),$$
so that 
\be \label{big}
 |\int xK_{i}\circ \f_{i}| \geq  \psi_{q}(t_{i}),
\ee

  Comparing the behavior of $\chi(t)$ and $\psi(t)$ in \eqref{kest1} and \eqref{cy21}, the estimates \eqref{small} and \eqref{big} give 
a contradiction for $t_{i}$ sufficiently large, i.e.~$\f_{i}$ sufficiently large in $\Co$. It follows that $\{\f_{i}\}$ is bounded in $\Co$, 
and the proof is then completed above. 

{\endproof}

\begin{remark}
{\rm The properness of the map $\pi_{0}$ corresponds to (and is in fact  equivalent to) general 
{\it a priori} estimates for solutions of \eqref{u}. Thus, as noted in the Introduction, the properness of $\pi_{0}$ 
implies that for any $K \in \cN\cap \cC$ with $dist_{C^{m-2,\a}}(K, \partial (\cN \cap \cC)) \geq \e > 0$, there is a constant 
$C = C(|K|_{C^{m-2,\a}}, \e, \a')$ such that, if $u$ is any solution of \eqref{u}, (i.e.~$\pi_{0}(u) = K$), then 
\be \label{apri}
|u|_{C^{m,\a'}} \leq C,
\ee
for any $\a' < \a$. 
}
\end{remark}

\begin{remark}
{\rm We note that the use of Astala's theorem in \cite{As} is crucial in the argument above. If $\eta_{i}$ are merely bounded 
in $W^{1,2}$, one only obtains an $L^{1}$ bound on $e^{2v_{i}}$ from \eqref{a4} and this is not enough to obtain control on 
$v_{i}$ from \eqref{vi}. In fact obtaining suitable control in this case relates to the exponential Sobolev embedding theorems 
of Trudinger-Moser-Aubin.  On the other hand, one only really needs a bound on $\eta_{i}$ in $W^{1,p}$ for some $p > 2$, 
not any $p < \infty$.  
}
\end{remark}

  The structure of the domain $\cK = \cN \cap \cC$ plays an important role in the further analysis of the curvature map 
$\pi$. Let $\{\cK^{i}\}$ denote the collection of path components of $\cK$, so that 
\be \label{ni}
\cK = \bigcup_{i}\cK^{i}.
\ee
As in \eqref{pii}, let $\cU^{i} = \pi^{-1}(\cK^{i})$. The map $\pi$ (or $\pi_{0}$) restricted to $\cU^{i}$ gives then a proper 
Fredholm map 
\be \label{pii2}
\pi^{i}: \cU^{i} \to \cK^{i},
\ee
into a connected target space. It is not easy to determine the number of components of $\cK$; we will see in the next 
section that $\cK$ has in fact infinitely many components.

   Since both $\cN$ and $\cC$ are open and dense in $C_{+}^{m-2,\a}$, the closure $\overline{\cK} = C_{+}^{m-2,\a}$. 
The point-set theoretic boundary 
\be \label{bound}
\partial \cK = C_{+}^{m-2,\a} \setminus \cK = \partial \cN \cup \partial \cC
\ee
consists of two main parts; $\partial \cN$ consists of functions $K$ such that $\D K(q) = 0$ for some critical point 
$q$ of $K$ with $K(q) > 0$, while $\partial \cC$ consists of functions $K$ for which there exist $p \in \{K = 0\}$ such 
that \eqref{criterion} holds at $p$ and $\D K(p) \leq 0$. Note that if $p$ is a non-degenerate critical point of $K \in \partial \cC$, 
then $p$ is necessarily a local maximum. Both of these boundaries will play an important role in the 
analysis to follow. 

\begin{remark}
{\rm  As noted above, while $\cN$ and $\partial \cN$ are invariant under the action of $\Co$, they are not invariant under 
the natural action $(\psi, K) \to K\circ \psi$ of $\Di^{m-2,\a}$ on $C^{m-2,\a}$. (On the other hand, $\cC$ and $\partial \cC$ 
are $\Di^{m-2,\a}$ invariant). It is easy to see that the sign of $\D K$ at a saddle point $q$ may be changed arbitrarily 
under the action of local diffeomorphisms acting only in a neighborhood of $q$. 

  Also, the number of critical points of $K$ (even Morse non-degenerate critical points) may change in a given path component 
$\cN^{i}$ of $\cN$. For instance, in a neighborhood of a local maximum of $K$, one may locally perturb $K$ by ``pushing down",  
creating two local maxima and a saddle point (mountain pass), while keeping $\D K < 0$ throughout the process. 
}
\end{remark}

  It is worthwhile to compare the space $\cK$ with the space $\cM$ of Morse functions in $C_{+}^{m-2,\a}$. 
The space $\cM$ is also open and dense in $C_{+}^{m-2,\a}$ and is $\Di^{m-2,\a}$ invariant. As in \eqref{bound}, 
one has $\partial \cM = C_{+}^{m-2,\a} \setminus \cM$ and there is a decomposition into $\Di^{m-2,\a}$ invariant 
path components 
\be \label{mi}
\cM = \bigcup_{j} \cM^{j}.
\ee
Now it is clear that there are infinitely many path components $\cM^{j}$, since for instance the number of 
non-degenerate critical points, which can be arbitrarily large, is constant on each component $\cM^{j}$. 
The boundary or wall $\partial \cM$, consisting of functions 
with at least one degenerate critical point, serves as a birth (or death) region for the creation of non-degenerate critical 
points. Under perturbation, a degenerate critical point in $\partial \cN$ can resolve for instance into either two critical points 
or zero critical points. For a detailed discussion of the decomposition \eqref{mi}, we refer to \cite{Ar}, \cite{Ni}. 

\medskip 

   We have the following structure of generic points in $\partial \cK$. 

\begin{proposition} The intersection $\cM \cap \partial \cK$ is open and dense in $\partial \cK$. Near a generic $K \in 
\cM \cap \partial \cK$, the boundary $\partial \cK$ is a smooth, codimension 1 hypersurface of $C_{+}^{m-2,\a}$, 
locally separating $C_{+}^{m-2,\a}$ into two components. In a neighborhood of any $K \in \cM \cap \partial \cK$, 
$\partial \cK$ is a region formed from the transverse intersection of a finite number of codimension 1 hypersurfaces. 
\end{proposition}

{\bf Proof:} It is obvious that $\cM \cap \partial \cK$ is open in $\partial \cK$. To prove the intersection is dense, 
we work first with $\partial \cN$. Let $K \in \partial \cN$ and let $K_{i}$ be a sequence of Morse functions with 
$K_{i} \to K \in C_{+}^{m-2,\a}$. If in a subsequence $K_{i} \in \partial \cN$, there is nothing more to prove, so suppose 
$i$ is large and $K_{i} \notin \partial \cN$. For each $i$ large, the function $K_{i}$ has a finite set $Q_{i}$ of non-degenerate 
critical points which are Hausdorff close to the critical locus $Z = \{q: |\nabla K|(q) + |\D K|(q) = 0\}$ of $K$. 
By means of small local perturbations by diffeomorphisms as discussed in Remark 2.12, one may perturb $K_{i}$ slightly, 
keeping $K_{i}$ Morse, and such that $|\nabla K_{i}| + |\D K_{i}| = 0$ on the saddle points of $Q_{i}$ while $|\Delta K_{i}| \leq 
\e_{i}$ on any local maxima and minima in $Q_{i}$; here $\e_{i} \to 0$ as $i \to \infty$ . Letting $i \to \infty$ then proves that 
$\cM \cap \partial \cN$ is dense in $\partial \cN$. The proof that $\cM \cap \partial \cC$ is dense in $\partial \cC$ is 
similar and straightforward. 

   Next, define 
$$F: C_{+}^{m-2,\a} \times S^{2} \to \bR^{3}$$
$$F(K, p) = (\nabla K(p), \D K(p)).$$
Thus $K \in \partial\cN$ if there exists $p$ such that $(K, p) \in F^{-1}(0)$. The map $F$ is smooth and the 
derivative in the first factor is $DF_{p}(\b) = \frac{d}{dt}F(K+t\b)(p)= (\nabla \b(p), \D \b(p))$, which is clearly 
surjective since $\b$ is arbitrary. Thus, by the implicit function theorem in Banach spaces, $\S = F^{-1}(0)$ is a 
codimension 3 smooth hypersurface of $C_{+}^{m-2,\a} \times S^{2}$. Of course $F^{-1}(0)$ may not be connected. 

  Now consider the projection $p: C_{+}^{m-2,\a}\times S^{2} \to C_{+}^{m-2,\a}$ of $\S$. Clearly 
\be \label{p}
p(\S) = \partial\cN.
\ee  
We claim that for $K \in \cM$, $\S$ is transverse to the vertical $S^{2}$ factors, in that  
\be \label{trans}
T\S \cap TS^{2} = 0.
\ee
To see this, one has $T\S = Ker DF$ and for $v \in T_{(K,p)}S^{2}$, $F_{*}(v) = \frac{d}{dt}(\nabla K(p+tv), \D K(p+tv)) = 
(D^{2}K(v), v(\D K))$. Thus $v \in T\S \cap TS^{2}$ if and only if $v \in Ker D^{2}K$ and $v \perp \nabla \D K$ at $p$. 
Since $p$ is a critical point of $K$ and $K$ is a Morse function, it follows that $v = 0$, which proves the claim. 

  Given any $K \in \cM \cap \partial\cN$, consider $\D K$ as a function on the set of critical points $p$ 
of $K$. Generically, a point $K \in \cM \cap \partial\cN$ has exactly one critical point $p$ where $\D K = 0$, 
with $\D K \neq 0$ at all other critical points. In this case, the transversality \eqref{trans} implies that $\partial \cN$ 
is a codimension one hypersurface near $K$. Such local hypersurfaces may have a finite number of intersections, 
corresponding to the vanishing of $\D K$ at multiple critical points. Since the critical points are isolated, it is 
clear that the intersection of $j$ such hypersurfaces is a submanifold of codimension $j$. 

  Similarly, a function $K \in \cM$ has a finite number of local maxima and by considering the behavior of $K$ on such 
local maxima, it is easy to see that $\cM \cap \partial \cC$ is the intersection of a finite number (generically 
one) of smooth hypersurfaces. 

{\endproof} 

  The proof of Proposition 2.10 shows that in general $\partial \cK$ has the structure of a rectifiable set of codimension 
one in $C_{+}^{m-2,\a}$. Namely, by definition, a set $Y \subset C_{+}^{m-2,\a}$ is rectifiable of codimension $s$ if $Y$ is 
contained in a countable union of sets of the form $p(G)$, where $p: C_{+}^{m-2,\a}\times \bR^{k} \to C_{+}^{m-2,\a}$ is 
projection on the first factor and $G$ is a codimension $s+k$ submanifold (a graph) in $C_{+}^{m-2,\a}\times \bR^{k}$, 
cf.~\cite{W1}, \cite{Be} for further details. Rectifiable sets of codimension $s$ are preserved under smooth Fredholm maps of 
index 0 (cf.~\cite{Be}). The characterization of $\partial \cN$ in \eqref{p} then shows that $\partial\cN$ 
is rectifiable of codimension one and similarly for $\partial \cC$.  

 Although Proposition 2.13 shows that  $\partial \cK$ {\it locally} separates $C_{+}^{m-2,\a}$, it is not easy to see to 
what extent $\partial \cK$ globally separates $\cK$ into distinct components. This will be discussed further in 
Sections 3 and 4. 

  Regarding $\cC$, it is clear that the space of positive curvature functions 
\be \label{cP}
\cP = \{K \in C^{m-2,\a}: K(x) > 0, \ \forall x \in S^{2}\},
\ee
is path connected and of course contained in $\cC$. 

\begin{proposition} The space $\cC \subset C_{+}^{m-2,\a}$ is path connected. 
\end{proposition}

{\bf Proof:}   Let $K(r)$, $r \in [-1,1]$ be a smooth path of functions in $C_{+}^{m-2,\a}$ connecting points $K(-1)$ 
and $K(+1)$ in $\cC$. Without loss of generality, assume $K(r)$ has only isolated critical points and assume $K(r)$ 
intersects $\partial \cC$ in the single point $K(0)$. Let $K_{\pm} = K(\pm \e)$ with $\e$ sufficiently small and assume 
$K_{\pm}$ are Morse functions. 

  We will use the following simple construction from Morse theory. Let $K$ be a Morse function. Then all pairs of local 
maxima of $K$ are connected by descending gradient flow lines to the collection of saddle points of $K$. Consider the 
collection $C$ of local maxima in the region $U_{0} = \{K < 0\} \subset S^{2}$. Each local maximum in $C$ is connected 
via a saddle point in $U_{0}$ to another local maximum in $C$ or to a local maximum (perhaps the absolute maximum) 
in the region $K > 0$. Here without loss of generality assume that $0$ is a regular value of $K$ so that the level set 
$L_{0} = \{K = 0\}$ consists of a finite collection of embedded simple curves, bounding the planar domain $U_{0} 
\subset S^{2}$. 

    Then by standard Morse theory  there is a smooth curve of functions $K_{t}$, leaving $K$ unchanged in the region 
$K \geq 0$, deforming $K = K_{0}$ to $K_{1}$ and for which $K_{1}$ has no local maxima in the region $K_{1} < 0$. 
Namely by pushing down along gradient flow lines, one may inductively cancel off any local maximum with a saddle 
point in $U_{0}$. (Of course this path passes from one Morse chamber to another, but this is irrelevant here). 
It follows that $K_{1}$ has only local minima and saddle points in the region $K_{1} < 0$. The path 
$K_{t}$ is contained in $\cC$, i.e.~it does not intersect $\partial \cC$. One may visualize this construction concretely 
by viewing $K$ as the height function $z$ of an embedding or immersion of $S^{2}$ into $\bR^{3}$; this is always 
possible by the main result in \cite{Ku}. 

  Now apply the construction above to the two points $K = K_{\pm}$. This gives a pair of paths joining $K_{\pm}$ to 
points $K_{\pm}'$ remaining within $\cC$. Next, combine these paths with the paths $K_{\pm}'(s) = K_{\pm}' + sc$ 
for $s \in [0,1]$ and $c > \max(|\min K_{+}'|, |\min K_{-}'|)$. The resulting paths join $K_{\pm}'$ to points $K_{\pm}'(1)$ 
in $\cP$ again remaining within $\cC$. The endpoints of these curves may then be connected within $\cP$. 

{\endproof}

Observe that any component $\cN^{i}$ of $\cN$ intersects $\cP$; namely as above, for any $K \in \cN^{i}$, choosing  
$c > |\min K|$, the path $K(s) = K + sc$, $s \in [0,1]$ connects $K$ within $\cN^{i}$ to a point in $\cP$. 

\begin{proposition}
Any path component $\cK^{i}$ of $\cK = \cN \cap \cC$ intersects $\cP$. Moreover, $\cK^{i}\cap \cP$ is path connected 
so that a path component $\cK^{i}$ of $\cK$ is uniquely determined by the path component $\cK^{i}\cap \cP$. 
\end{proposition}

{\bf Proof:} The proof is a refinement of the proof of Proposition 2.14. Choose any Morse function $K \in \cM \cap \cK$. 
First note that if $K$ has no negative local maxima, i.e.~no local maxima $p$ with $K(p) < 0$, then adding a curve of 
positive constants to $K$ as above connects $K$ to $\cP$ within $\cK$, i.e.~within a path component of $\cK$. 

  Next suppose $p$ is a negative local maximum of $K$ which is connected via a gradient flow line to a negative saddle 
point $q$, in the sense that $\D K(q) < 0$; necessarily $K(q) < 0$. Then as in the proof of Proposition 2.14, one can 
deform $K = K_{0} \in \cK$ to $K_{1} \in \cK$ along a curve by pushing down along gradient flow lines to cancel $p$ and $q$, 
while maintaining $\D K < 0$ at the curve of saddle points. As above this requires crossing a wall between Morse chambers 
so that a saddle point becomes Morse degenerate, i.e.~$det D^{2}K(q_{t}) = 0$. However, one may maintain $\D K(q_{t}) < 0$ 
on the curve $q_{t}$ of saddle points. Such curves thus remain in $\cK$. One may iterate this process inductively to cancel off all 
such pairs. This connects $K$ to a point $K'$ to which one may apply the argument above.  

  Thus suppose $K$ has a negative local maximum $p$ without negative (i.e.~$\D K < 0$) saddle points in $U_{0}$ connected 
to $p$. In this case, we first create a new pair of positive local maxima with a new negative saddle point starting from $K$. 
Namely the given function $K$ has a positive absolute maximum say at $p_{0}$. By locally pushing down near $p_{0}$ 
as in Remark 2.12, one may create within $\cK$ a new function $K'$ with two positive maxima and a saddle $q'$ connecting 
them, with $\D K'(q') < 0$ with $K = K'$ on $U_{0}$, were $U_{0} = S^{2}\setminus V_{0}$ and $V_{0}$ is a small 
neighborhood of $p_{0}$. The value $K'(q')$ of the saddle $q'$ may then be pushed down further so that $K'(q') < K(q)$, 
for any saddle $q$ of $K$ in the region $U_{0}$, again remaining in $\cK$. Replacing $K$ by $K'$, one may then apply 
the construction in the previous paragraph to cancel a negative local maximum of $K' = K$ with the negative saddle point 
$q'$. Again, this construction may be iterated inductively to cancel all negative local maxima of while remaining within 
$\cK$. In other words, $K$ is connected by a curve in $\cK$ to another Morse function $\hat K$ without negative local 
maxima. One may then connect $\hat K$ to $\cP$ by pushing up by positive constants as above. 

  The second statement follows from the stability or continuity of the construction above. Thus, suppose $K_{1}$ and $K_{2}$ 
are two points in $\cK^{i}\cap \cP$ connected by a smooth path $K(s)$ in $\cK^{i} \subset \cN \cap \cC$. The process 
above applied at $K_{1}$ may be applied continuously along $K(s)$ to deform the path $K(s)$ into $\cK^{i}\cap \cP$.  

{\endproof}

\section{degree computations} 

  The purpose of this section is to prove Theorems 1.2 and 1.3. 
  
\medskip   
  
 The map $\pi_{0}$ in \eqref{pi0} or \eqref{pi01} is a proper Fredholm map of index 0 between open subsets of Banach spaces. 
Proper Fredholm maps $F: X \to Y$ of index 0 between oriented Banach spaces or oriented Banach manifolds 
$X$, $Y$ with $Y$ path connected have a $\bZ$-valued degree defined as 
$$deg \, F = \sum_{x\in F^{-1}(y)} sign(D_{x}F) \in \bZ,$$
wbere $y$ is any regular value of $F$ and the sign of $DF(x)$ is $\pm 1$ according to whether $DF(x)$ preserves 
or reverses orientation at $x$, cf.~\cite{ET} or more recently, \cite{BF1}, \cite{BF2} for example. A point $y \in Y$ is a 
regular value of $F$ if $D_{x}F$ is an isomorphism, for all $x \in F^{-1}(y)$. Note that the mod 
$2$ reduction of $deg_{\bZ}F$ is the (unoriented) Smale degree \cite{Sm}. The domain $\cU$ of $\pi_{0}$ is 
an open subset of a Banach space, and so is clearly orientable, as is the target space $\cN \cap \cC \subset C_{+}^{m-2,\a}$. 

\medskip 

{\bf Proof of Theorem 1.2.}

\medskip 

  To calculate the degree of 
$$\pi^{i}: \cU^{i} \to \cK^{i},$$
as in \eqref{pii2}, recall that Proposition 2.15 shows that each component $\cK^{i}$ of $\cK$ intersects $\cP$ and $\cK^{i}$ itself 
is path connected. The degree of $\pi^{i}$ is then the same as the degree of the restricted map 
$$\w \pi^{i}: \cU_{\cP}^{i} \to \cK^{i} \cap \cP,$$
where $\cP$ is the space of positive curvature functions as in \eqref{cP} and $\cU_{\cP}^{i} = (\pi^{i})^{-1}(\cK^{i}\cap \cP)$. 

 The degree of $\w \pi^{i}$ is calculated by Chang-Gursky-Yang in \cite{CGY}, by an elegant argument based on the 
 Poincar\'e -Hopf theorem and is given by the formula \eqref{degform0}. Hence \eqref{degform0} holds for $\pi^{i}$.   

{\endproof}

\begin{remark}
{\rm Note  that the formula for $deg \, \pi^{i}$ given by \eqref{degform0} changes, generically by one, when passing through 
the boundary region $\partial \cC \cap \cN$, i.e.~on a curve $K(t)$ on which a local maximum of $K(t)$ passes through 
zero. Similarly it changes generically by one on curves $K(t)$ passing through the boundary region $\partial \cN \cap \cC$, 
so that $\D K(t)$ changes sign at a saddle point of $K(t)$. Note that the sign of $\D K$ can be changed or altered just 
by means of local diffeomorphisms of $S^{2}$ supported near the saddle points, as in Remark 2.12.
}
\end{remark}
 
   We note that the results above in Sections 2 and 3 do not use any variational formulation of the curvature equation 
\eqref{u}. We now turn to the proof of Theorem 1.3 which gives another calculation of $deg \, \pi^{i}$ using the 
functional $J$ from \eqref{J}. 

   Thus consider 
\be \label{I2}
J: C^{m,\a}\times C_{+}^{m-2,\a} \to \bR,
\ee
$$J(u, K) = \int_{S^{2}}(|d u|^{2} + 2u)dv_{0} - \log (\int_{S^{2}}Ke^{2u}dv_{0}).$$
Note that $J$ is invariant under the translations $u \to u + c$. The variation $D_{1}J$ of $J$ at $(u, K)$ with respect 
to the first variable is given by 
\be \label{1var}
{\tfrac{1}{2}}D_{1}J =  -\D u  + 1 - \k Ke^{2u},
\ee
where $\k = (\int Ke^{2u}dv_{0})^{-1}$. Hence, critical points of $J$ with respect to variations of $u$, subject to the 
constraint $\k = 1$, are exactly the solutions of \eqref{u}. It is worth noting that the structure, for example the topology, 
of the constraint manifold 
$$\cC_{K} = \{u \in C^{m,\a}: \int_{S^{2}}Ke^{2u}dv_{0} = 1\},$$
depends on $K \in C_{+}^{m-2,\a}$. Of course there are no critical points with respect to 
variation in the second or $K$-variable. By Theorem 1.1, the set of critical points of $J_{K}$ (satisfying the constraint 
$\k = 1$) is compact for any fixed $K \in \cK$. If $K$ is a regular value of $\pi_{0}$, there are only finitely many critical 
points of $J_{K}$ and all are non-degenerate, (cf.~Proposition 5.3 below). 

  Given a critical point $u$ of $J_{K}$, the index $ind_{u}$ of $J_{K}$ at $u$ is the maximal dimension on which the second 
variation $D^{2}J_{K}$ with respect to $u$ is negative definite. It is easy to see that $ind_{u}$ is finite, for any fixed $u$. 
We note that it is well-known \cite{Ho} the infimum of $J_{K}$ is never achieved unless $K = const$. 

\medskip 

{\bf Proof of Theorem 1.3.} 

\medskip 

This result is by now an essentially standard result for variational problems satisfying compactness 
conditions as in Theorem 1.1. In the following, we adapt an elegant approach of White \cite{W1}, \cite{W2} to this 
setting. 

  The right-hand side of \eqref{deg2} is defined for each regular value $K$ and the main issue is to prove the expression is 
well-defined, i.e.~independent of the choice of $K \in \cK^{i}$. Thus let $K(0)$ and $K(1)$ be any two regular values of $\pi^{i}$ 
and let $K(t)$, $t \in [0,1]$, be an oriented smooth curve in $\cK^{i}$ joining them. Let $u(t) = (\pi_{0})^{-1}(K(t))$. 
As shown in Theorem 5.5 below (which is independent of the considerations here), the set of singular points $\S \subset \cU^{i}$ 
is a stratified space of codimension at least one. Hence, if necessary, $K(t)$ may be perturbed slightly so that $K(t)$ is transverse 
to $\pi^{i}$ and thus $u(t)$ is a collection of smooth disjoint curves in $\cU^{i}$ connecting regular points over $K(0)$ and 
$K(1)$, cf.~also \cite{Sm}. Define an orientation on $u(t)$ by defining $\pi_{0}$ to be orientation preserving near any regular 
point $u(t)$ if the index $ind_{u(t)}$ is even, orientation reversing if the index is odd. Given this orientation on curves, the 
mapping degree of 
$$\pi_{0}: u(t) \to K(t),$$
between these 1-dimensional manifolds is then given exactly by \eqref{deg2}. It suffices then to show that this definition of 
orientation is well-defined, i.e.~consistent, in that when $u(t)$ passes through a critical point of $\pi_{0}$, the local orientation 
of $\pi_{0}$ changes according to the change of the index. In the following, we assume without loss of generality that $u(t)$ is a 
connected curve. 

Suppose $u(t_{0})$,  $t_{0} \in (0,1)$ is a critical point of $\pi_{0}$, so that the Hessian $D^{2}J_{K}$ at $u(t_{0})$ has 
a non-trivial null-space $N$, i.e.~the nullity $dim N$ of $J_{K}$ at $u(t_{0})$ is non-zero, cf.~also Proposition 5.3 below.  
Without loss of generality, assume $dim N = 1$; the arguments below apply to each line in $N$. The following 
computations take place for $t$ near $t_{0}$. 

  It is trivial to see that the map $F: C^{m,\a}\times C_{+}^{m-2,\a} \to C^{m-2,\a}$ given by $F(u, K) = \frac{1}{2}D_{1}J = 
1 - \D u - \k Ke^{2u}$ is a submersion. Hence $F^{-1}(0)$ is a closed submanifold; clearly $F^{-1}(0) = {\rm graph} \,\pi$. 
The submanifold $F^{-1}(N)$ is a 1-dimensional thickening of $F^{-1}(0)$ and one may thicken the curve 
$(u(t), K(t))$ to a 2-parameter family 
$$\s (t,s) = (u(t,s), K(t,s)) \in F^{-1}(N)  \ \ {\rm with} \ \ K(t,s) = K(t).$$ 
By construction, one has 
\be \label{intf}
\frac{d}{ds}J(\s(t,s)) = \int_{S^{2}}\<F(\s(t,s)), \frac{du}{ds}\>dv_{+1}.
\ee
Clearly $\frac{d}{ds}J(\s(t,s))|_{s=0} = 0$. By construction $F(\s(t,s)) \in N$ and $\frac{du}{ds} \in N$ are non-zero 
for $s \neq 0$. In the $(t, s)$ plane $P$, the locus where $\partial J(t,s)/\partial s = 0$ (with $J(t,s) := J(\s(t,s))$) is given locally 
by the graph $s = u(t)$ of the curve $u(t)$. In particular, $(t, u(t))$ is locally the boundary of the open domain 
$U = \{\partial J /\partial s > 0\}$ in $P$: 
$$u(t) = \partial U.$$
If $ind_{u_{t}}$ is even for $t$ near $t_{0}$, $t \neq t_{0}$, $u(t)$ is given the boundary orientation induced by $U$ while if 
$ind_{u(t)}$ is odd near $t_{0}$, give $u(t)$ the opposite orientation. It is now simple to see that this choice agrees with the 
definition above. Namely, the point $t_{0}$ is a critical point for the map $\pi_{2}\circ u: I \to \bR$, where $\pi_{2}$ is projection 
onto the $s$-factor. If this critical point is a folding singularity for $\pi_{2}\circ u$, then the index of $u(t)$ changes by 1 in 
passing through $u(t_{0})$ and reverses the orientation of $\pi_{2}(u(t))$, (exactly as in the case of the standard folding 
singularity $x \to x^{2}$). If on the other hand $\pi_{2}\circ u$ does not fold with respect to $\pi_{2}$ (so that one has an 
inflection point) then the index of $u(t)$ does not change through $t_{0}$ and $\pi_{2}$ maps $u(t)$ to $\pi_{2}\circ u(t)$ in an 
orientation preserving way. 

  This shows that the degree \eqref{deg2} is well-defined. The fact that it equals the degree in \eqref{deg2} up to an overall sign 
follows from the basic uniqueness properties of such $\bZ$-valued mappings, cf.~\cite{BF1}. This completes the proof. 
 
{\endproof}

\section{Discussion on the Image of $\pi$.}

  In this section we summarize some of the main results of the previous sections and prove several further results on the 
behavior of $\pi$ near the boundary $\partial \cK$. We also discuss the structure of the components of $\cK$. These 
issues provide a bridge to issues discussed in the next section. 

\medskip 

  It follows from the results in Sections 2 and 3 that the boundary $\partial \cK$ is a closed rectifiable set of 
codimension 1, separating $C_{+}^{m-2,\a}$ into an infinite collection of components $\cK^{i}$ on which the proper maps 
$$\pi^{i}: \cU^{i} \to \cK^{i}$$
have a well-defined degree. For any $K$ in a component $\cK^{i}$ of non-zero degree, \eqref{u} is solvable (since 
$\pi^{i}$ is surjective), and there is an effective, signed, count on the number of solutions for generic $K$. There 
are no restrictions on the structure of the set of critical points for $K$ in such components $\cK^{i}$; for example, 
critical points of $K$ need not be isolated. 

   The only possible regions where \eqref{u} may not have a solution are the components of degree zero and 
the boundary region $\partial \cK = \partial\cN \cup \partial\cC$. 

\medskip 

   First we discuss some general structural results for $\cK$. 

\begin{lemma}
The degree of $\pi_{0}$ may take on any value in $\bZ$, i.e.~for all $n \in \bZ$, there exists a component $\cK^{i}$ 
of $\cK$ such that 
$$deg \, \pi^{i} = n.$$
\end{lemma}

{\bf Proof:} Let $K$ be a Morse function in $\cK$. Let $M$ and $m$ denote the number of local maxima and minima 
respectively, and let $s^{+}$, $s^{-}$ denote the number of saddle points with $\D K > 0$ and $\D K < 0$ respectively. 
Then the Morse Lemma gives $M - (s^{+} + s^{-}) + m = \chi(S^{2}) = 2$, while the degree formula \eqref{degform0} 
gives $deg \, \pi_{0} = M - s^{-} - 1$. 

  If all saddles have $\D K > 0$, so $s^{-} = 0$, then $deg \, \pi_{0} = M - 1$, while if all saddles have $\D K < 0$, then 
$deg \, \pi_{0} = M - (s^{+} + s^{-}) - 1 = 1 - m$. It is clear that for any $M$ there is a Morse function with $M$ local maxima and 
similarly for any $m$, there is a Morse function with $m$ local minima.  

  Finally, as in Remark 3.1, at any given saddle point $q$, by using local diffeomorphisms supported near $q$, one 
may change the sign of $\D K$ arbitrarily. It follows that $deg \, \pi_{0}$ takes on all values in $\bZ$. 

{\endproof}

\begin{remark}
{\rm It is unknown whether the set of components $\cK^{i}$ of a fixed degree is connected, (although this seems unlikely). 
In view of Proposition 2.15, consider the collection of components $\cK^{i}\cap \cP$, $i$ varying, of a fixed degree, say 
$n = M - s^{-} - 1$ and let $K_{i} = K$ be a Morse function in $\cK^{i}\cap \cP$. Each such $K$ has a finite number of 
critical points and the flow lines of $\nabla K$ flow between these critical points. Any pair of local maximum points $p_{1}$, 
$p_{2}$ of $K$ connected to a saddle point $q$ of $K$ with $\D K(q) < 0$ can be canceled by a path $K(t)$, $t \in [0,1]$, of 
Morse functions staying in the given component $\cK^{i}\cap \cP$, leaving a single maximum point. Thus both $M$ and $s^{-}$ can 
be reduced by one along such paths in $\cK^{i}\cap \cP$. In the same way, any pair of local minima of $K$ can be canceled within 
$\cK^{i}\cap \cP$ with a saddle point with $\D K > 0$ leaving a single minimum point. Thus, within each component $\cK^{i}\cap \cP$, 
one may find Morse functions $K_{i}' \in \cK^{i}\cap \cP$  with the minimal numbers $M_{min}$, $m_{min}$ of local maxima and minima 
respectively - among Morse functions in $\cK^{i}\cap \cP$. The numbers $M_{min}$ and $m_{min}$, and so also the minimal number of 
saddle points $s_{min}$, are thus the same for the collection of components $\cK^{i}\cap \cP$. Note also that the space of Morse functions 
with a fixed number of local maxima and minima is path-connected in $\cM$, cf.~\cite{Ku} for example. However it is not clear that 
there must exist paths connecting such points in $\cK^{i}\cap \cP$, i.e.~paths on which the sign of $\D K$ at the saddle points 
remains fixed. 

  On the other hand, note that the argument above proves that there is only one component $\cK^{0}$ of degree zero, since 
$M_{min} = m_{min} = 1$, $s_{min} = 0$ on $\cK^{0}$. 

}
\end{remark}
 
   Next we show that there are solutions of \eqref{u} in each of the regions $\cK^{0}$ (the degree zero component) 
and $\partial \cK$ mentioned above. 
   
\begin{proposition} 
The domain $\cK^{0}$ is non-empty, i.e.~there exist $u$ with $K_{u} = \pi(u) \in \cK^{0}$.
\end{proposition}

{\bf Proof:} It is well-known, cf.~\cite{KW1}, that given any function $K \in C_{+}^{\infty}$ there is a smooth metric $g$ on $S^{2}$ 
such that $K_{g} = K$. By the uniformization theorem $g = \psi^{*}(e^{2u}g_{+1})$, for some diffeomorphism $\psi$ and function 
$u$. The metric $\w g = e^{2u}g_{+1}$ is conformal to $g_{+1}$ with Gauss curvature $K_{\w g} = K\circ \psi^{-1}$. 

  Now choose $K$ to be any Morse function in $\cK^{0}$ with exactly two critical points, necessarily a maximum and a 
minimum. This property is invariant under diffeomorphisms, and necessarily $\w K = K_{\w g} \in \cK^{0}$. Then $\pi(u) = \w K$. 

{\endproof} 

  Although Remark 4.2 shows that $\cK^{0}$ is only one component among infinitely many others on which $\pi_{0}$ is onto, 
the next result shows that, in a certain sense, the degree zero component $\cK^{0}$ is the largest. 

\begin{proposition}
For any given $K$, the set of linear functions $\ell$ such that
\be \label{largel}
K + \ell \in Im \, \pi
\ee
is a bounded set in $\bR^{3}$. Hence, for any given $K$, for any $\ell$ sufficiently large, $K + \ell \in \cK^{0}$ and $K + \ell$ 
is not solvable for $u$ in \eqref{u}. 
\end{proposition}

{\bf Proof:} The proof is (naturally) based on the Kazdan-Warner obstruction \eqref{kw}. First if $K = const$, then 
$K + \ell \notin Im \pi$ for any $\ell \neq 0$ by \eqref{kw}. Thus we assume $K \neq const$. Any given $\ell$ may 
be written in the form $\ell = c\<p, \cdot\>$, where $p \in S^{2}$, $c > 0$ and $\< \cdot, \cdot\>$ is the Euclidean 
dot product. As noted following \eqref{comm1}, $K \in Im \, \pi$ if and only if $K_{\f} = K\circ \f \in Im \, \pi$ for any 
$\f \in \Co$. Given $p$ as above, choose $\f$ such that $\min K_{\f}$ is achieved at $-p$ while $\max K_{\f}$ is achieved 
at $p$; this is always possible since $K \neq const$.  Now choose the conformal vector field $X = \nabla \ell_{0}$ and 
let $\ell_{\f} = \ell \circ \f$. 

   Then if $K_{\f} + \ell_{\f} \in Im \, \pi$, one must have  
$$\int_{S^{2}}X(K_{\f} + \ell_{\f})dv_{g} = \int_{S^{2}}[X(K_{\f}) + cX(\ell_{0}\circ \f))]dv_{g} = 0,$$
for some volume form $dv_{g}$. Note that $X(\ell_{0}\circ \f) = \<D\f(\nabla \ell_{0}), X\>$ with respect to 
$g_{+1}$. Since $\f$ is conformal, $X(\ell_{0}\circ \f)_{x} = \chi(x)X(\ell_{0})_{\f(x)}$ where $\chi = \sqrt{det D\f} 
> 0$. One has $X(\ell_{0}) = \xi = \sin^{2}r \geq 0$. This gives 
$$\int_{S^{2}}X(K_{\f} + \ell_{\f})dv_{g} = \int_{S^{2}}[X(K_{\f})(x) + c\chi(x) \xi(\f(x))]dv_{g} = 0.$$
Using the change of variables formula, this may be rewritten as 
\be \label{ellcont2}
\int_{S^{2}}X(K_{\f} + \ell_{\f})dv_{g} = \int_{S^{2}}[X(K_{\f})(x) + c\w \chi(x) \xi(x)]dv_{g} = 0.
\ee
where $\w \chi = \chi \circ \f^{-1}det D(\f^{-1})$is a positive functions on $S^{2}$. One has 
$0 \leq \xi\leq 1$, with $\xi$ vanishing only at the poles $\pm p$. If $c$ is sufficiently large, then $X(K_{\f}) + 
c\w \chi \xi > 0$ away from the poles $\pm p$, since the term $X(K_{\f})$ is fixed. On the 
other hand, at and sufficiently near the poles, $X(K_{\f}) \geq 0$ by construction and so again $X(K_{\f}) + 
c\w \chi\xi \geq 0$. This contradicts \eqref{ellcont2} (for any area form $dv_{g}$) if $c$ 
is sufficiently large.  
 
{\endproof} 

   Next, even though $\pi_{0}$ is not necessarily proper near $\partial \cK$, we show there are solutions in both parts of 
$\partial \cK$. 

\begin{proposition} 
The domain $\partial_{1}\cU : = \pi^{-1}(\partial \cN)$ is non-empty, i.e. ~there exist $u$ with $K_{u} = \pi(u) \in 
\partial \cN \subset \partial \cK$.
\end{proposition} 

{\bf Proof:}   There may be many proofs of this result, but we work here with a particular and simple, natural 
construction of a large class of functions $K \in Im \, \pi \cap \partial \cN$. 

  Namely, consider the behavior of $K$ on the space of $u$ which are eigenfunctions of the Laplacian, so that 
$\D u = -\l u$, $\l > 0$. These are the spherical harmonics, or equivalently, restrictions of harmonic polynomials 
in $x$, $y$, $z$ on $\bR^{3}$ to $S^{2}$. From \eqref{u}, for such $u$ one has 
\be \label{Keig}
K = e^{-2u}(1 + \l u).
\ee
Hence 
\be \label{Kcrit}
\nabla K = e^{-2u}((\l - 2) - 2\l u)\nabla u,
\ee
so that the locus of critical points of $K$ is the critical locus of $u$, together with the level set $u = (\l - 2)/2\l$. 
At the critical locus, one has 
\be \label{Klap}
\D K = - e^{-2u}((\l - 2) - 2\l u)\l u.
\ee
For any harmonic homogeneous polynomial $u$ of degree at least 2, it is easy to see that $0$ is always a 
critical value of $u$. Hence, for $K$ generated from arbitrary eigenfunctions of the Laplacian, one always has 
$$K \in \partial \cN.$$ 
Thus if $E_{\l}$ denotes the $\l$-eigenspace of the Laplacian on $S^{2}(1)$, then for any $\l$, 
$$\pi(E_{\l}) \subset \partial \cN.$$
Note that \eqref{Keig} and \eqref{Klap} imply that 
$$\pi(E_{\l}) \subset \cC.$$
{\endproof}

\begin{remark}
{\rm If $\l > 2$, one may restrict $u$ to the open ball in the eigenspace $E_{\l}$ where $\max u < (\l-2)/2\l$. For such $u$, 
the critical points of $K$ are exactly the critical points of $u$. 

  It is worth discussing a simple, specific example in more detail. Thus let $u = \frac{1}{4}(x^{2} - y^{2})$. This is an 
eigenfunction with $\l = 6$; note that $\frac{1}{4} < \frac{1}{3} = \frac{\l-2}{2\l}$. The critical points of $u$, 
equal to the critical points of $K$, come in three pairs. 

   The first pair is $x = \pm 1$, $y = z = 0$, giving $u = \frac{1}{4}$ and $\D K < 0$. This gives a pair of maxima 
of $u$ and $K$. The second pair is given by $y = \pm 1$, $x = z = 0$, giving $u = -\frac{1}{4}$ and a pair of minima 
of $u$ and $K$ with $\D K > 0$. The third pair $z = \pm 1$, $x = y = 0$ has $u = 0$ and gives a pair of saddle points of 
$u$ and $K$ with $\D K = 0$. A simple computation shows that $D^{2}K$ is non-degenerate at $z = \pm 1$, while 
$D^{2}K$ is degenerate at the other critical points. 

  Now consider perturbations of $u$ inducing perturbations of $K$ into $\cN$. Let $v = \eta w$ where $\eta$ is a cutoff 
function supported near $z = \pm 1$ with $\eta \equiv 1$ in a small neighborhood $U$ of the poles, and $w$ is another 
eigenfunction of $\D$ with eigenvalue $\l =  6$, say $w = z^{2} - y^{2}$. The computation of $K$ remains the same for 
$u + \e v$ in $U$ and the poles $z = \pm 1$ are still critical points.  By choosing $\e > 0$ or $\e < 0$, the poles resolve 
independently into non-degenerate saddles with $\D K \neq 0$; one may arrange $\D K > 0$ or $\D K < 0$ at a given 
pole by choosing $\e > 0$ or $\e < 0$ suitably. Outside of $U$, the $\e v$ perturbation does not change the structure 
of the other critical points of $K$. 

  This shows that for suitable perturbations of $u$, the curvature $K \in \partial \cN$ may be perturbed into a degree 0 
component of $\pi_{0}$,  (choosing one perturbation of the saddle to $\D K < 0$ and the other to $\D K > 0$), or into 
a degree 1 component, (choosing both saddle perturbations to $\D K > 0$), or into a degree $-1$ component, (choosing  
both saddle perturbations to $\D K < 0$). 

  Although $K$ here is not a Morse function, consider instead the eigenfunctions $u = x^{2} - (1+\e)y^{2} + \e z^{2}$. 
The same considerations apply as above and now $K = K(u)$ is a Morse function for $\e \in (0,1)$, so $K \in 
\cM \cap \partial \cN$. In connection with Proposition 2.13, such $K$ lie in the intersection of two regular local 
hypersurfaces of $\partial \cN$ and the structure of $\cN$ near $K$ is described by the wall-crossing above.  

   A similar structure holds for general eigenfunctions of the Laplacian, for any $\l$. 
}
\end{remark}

   A stronger existence result holds in the region $\partial \cC \subset \partial \cK$. The following result corresponds to 
Theorem 1.5. Write 
$$\partial \cC = \cup_{n,m}\partial \cC_{m}^{n},$$
where $\partial \cC_{m}^{n} = \partial \cC \cap \bar \cK^{m} \cap \bar \cK^{n}$ with $m < n$ and $\cK^{m}$, $\cK^{n}$ are 
components of degree $m$, $n$ respectively. Thus a curve $K(t)$ passing from $\cK^{m}\cap \cN$ 
to $\cK^{n}\cap \cN$ passes through $\partial \cC_{m}^{n}$, increasing the degree by $n - m$. 
 
\begin{proposition}
The map $\pi$ is onto $\partial \cC_{m}^{n}\cap \cM$ for any $(m, n)$ with $m \neq 0$, $m < n$.
\end{proposition}

{\bf Proof:} As discussed in Remark 2.7, the map $\pi$ is not proper at $\partial \cC$, due to the formation of (unbounded) 
bubble solutions as a local maximum of $K$ passes from small positive values to small negative values. This corresponds to 
a decrease in the degree of $\pi^{0}$ by one. 

   However, such bubbles do {\it not} form when passing in the reverse direction, from small negative to small positive values 
where the degree increases by one, at least within the space of Morse functions $K \in \cM$. To prove this, suppose $K \in \partial \cC 
\cap \cM$ so that there exists $x$ with $K(x) = 0$ and $x$ is a local maximum of $K$. Thus there exists $\d > 0$ such that 
$K(y) \leq K(x) = 0$ for all $y \in D_{x}(\d)$ where $D = D_{x}(\d)$ is the disc of radius $\d$ in $(S^{2}(1), g_{+1})$. Without 
loss of generality, assume there exists $\e_{0} > 0$ such that $K \leq -\e_{0} < 0$ on $\partial D$. Suppose there is a 
sequence $K_{i} \in \cC\cap \cN$ with $K_{i} = \pi(u_{i}) \to K$ in $C_{+}^{m-2,\a}$ such that $x$ is a local maximum 
of $K_{i}$ with $K_{i}(x) = -\e_{i} \to 0$ and $K_{i}(y) \leq 0$ for $y \in D_{x}(\d)$. We claim that $area_{g_{i}}D$ remains 
bounded as $i \to \infty$. Namely, since the curvature $K_{i}$ is non-positive on the simply connected disc $D$, the 
geodesic circles $(S_{i}(r), g_{i})$ of radius $r$ about $x$ with respect to $g_{i}$ are convex. Letting $\ell_{i}(r)$ denote 
the $g_{i}$-length of $S_{i}(r)$, it follows that $\ell_{i}'' \geq 0$ where the derivative is with respect to geodesic 
$g_{i}$-distance $r$. If $area_{g_{i}}D \to \infty$, it follows that the distance function $r_{i}$ becomes 
unbounded on $D$ as $i \to \infty$ and hence also 
\be \label{ellinf}
\ell_{i} \to \infty \ \ {\rm on} \ \ D. 
\ee
Now by \eqref{abound1.5}, the area of $g_{i}$ is uniformly bounded near the 
boundary $\partial D$ and in particular the $g_{i}$-length of $S = \partial D$ is uniformly bounded. This contradicts 
\eqref{ellinf}, which proves the existence of a bound on $area_{g_{i}}D$. It then follows from Theorem 2.3 and 
the proof of Theorem 2.9 that the sequence $g_{i}$ is bounded and, in a subsequence, $g_{i} \to g$, $u_{i} \to u$, in 
$C^{m,\a}$. This proves the claim above. 

  Now if $K_{i} \in \cK^{m}$ with $m \neq 0$, then since $\pi^{m}$ is onto, $K_{i} = \pi(u_{i})$, for some sequence 
$u_{i} \in \cU^{m}$. For $K = \lim K_{i} \in \partial \cC$ as above, suppose that $K \in \partial \cC \cap \bar \cK^{n}$, 
for $n > m$, so that $K \in \partial \cC \cap \bar \cK^{m}\cap \bar \cK^{n}$. This corresponds to an increase of the 
degree and the argument above proves that $K \in Im \, \pi$. 

{\endproof}

\section{Domain and Range Structure} 

  In this section, we discuss a number of aspects of the structure of the set of regular and singular points of $\pi$. 
The discussion is independent of the global results discussed in previous sections, until the analysis beginning with 
and following Proposition 5.7, which discuss global issues again. 

\medskip 
   
  Given a Fredholm map $F: X \to Y$ of index 0 between Banach manifolds, recall that $x$ is a regular point of 
$F$ if the linearization $D_{x}F$ is an isomorphism. The inverse function theorem implies that $F$ is a 
local diffeomorphism in a neighborhood of a regular point. A point $y \in Y$ is a regular value if all points 
$x \in F^{-1}(y)$ are regular points. By definition, any $y \notin Im F$ is a regular value. 
A point $x \in X$ is a singular point if it is not a regular point. 

  Let 
$$\S \subset C^{m,\a}$$
be the closed set of singular points of $\pi$, and $\cO$ the open set of regular points of $\pi$, so that 
$$C^{m,\a} = \cO \cup \S.$$ 
If $u$ is a singular point of $\pi$, so that $Ker D_{u}\pi \neq 0$, then there exists a deformation $u'$ of $u$ 
which preserves the Gauss curvature to first order: $K' = \frac{d}{dt}K(u + tu')_{t=0} = 0$. Recall the basic 
defining equation \eqref{u}, i.e. 
$$e^{2u}K = 1 - \D u.$$
Differentiating in the direction $u_{t} = u + tu'$ and setting $v =u'$, one has $e^{2u}(K' + 2vK) = -\D v$. 
This gives: 

\begin{lemma} 
A point $u \in C^{m,\a}$ is a critical point of $\pi$ if and only if there exists a variation $u' = \o$ such that 
\be \label{2}
\D \o + 2Ke^{2u}\o = 0. 
\ee
\end{lemma}

{\endproof}

  It is worth noting that by \eqref{comm1}, the sets of regular and singular points, as well as regular and singular values, 
are invariant under the action of $\Co$. 

\begin{remark}
{\rm  Since $\D_{\mu^{2}g} = \mu^{-2}\D_{g}$ in dimension two, one easily finds  
\be \label{dpi}
D_{u}\pi (v) = -\D_{g}v - 2Kv,
\ee
for $g = e^{2u}g_{+1}$. In case $K > 0$, this can be rewritten as 
\be \label{dpi2}
D_{u}\pi = -K(\D_{\w g}v + 2v),
\ee
where $\w g = Kg$. Thus $u$ is a critical point of $\pi$ if and only if $2$ is an eigenvalue of the 
conformal Laplacian $\D_{\w g}$. By a well-known theorem of Hersch \cite{H}, for any metric $\w \g$ on $S^{2}$ 
with $area(S^{2}, \w \g) = 4\pi$, one has 
$$\l_{1} \leq 2,$$
with equality if and only if $\w \g$ is of constant curvature $1$. Thus the conformal metrics $\w g$ in \eqref{dpi2} 
have $2$ as a $2^{\rm nd}$ (or higher) eigenvalue of the Laplacian except when $g$ has constant curvature. 

}
\end{remark}     
   
\begin{proposition}
The function $K$ is a regular value of $\pi$ if and only if all critical points of the functional $J_{K}$ in \eqref{J} 
are non-degenerate. 
\end{proposition}

{\bf Proof:} If $u$ is a critical point of $J_{K}$ giving a solution of \eqref{u}, then \eqref{1var} implies that 
$D^{2}J_{K}(v) = -\D v - 2Ke^{2u}v = e^{2u}(-e^{-2u}\D v - 2Kv) = e^{2u}K'(v)$. Thus, $v \in Ker D^{2}J_{K}$ if 
and only if $DK(v) = 0$, which proves the result.   

{\endproof}

  At a critical point $u$ of $\pi$, let $N$ denote the solution space of \eqref{2} and let $d = d_{u} = dim N$. 
The space $N = N_{u}$ equals the null-space of $D^{2}J_{K}$ at $u$, and $d$ is the nullity of the critical point $u$. 
By \eqref{dpi}, $Im D_{u}\pi$ consists of functions of the form $\D_{g}v + 2Kv$ where $g$ and $K$ are fixed and 
$v$ varies over $C^{m,\a}$. In particular 
\be \label{perp}
N  \perp Im D_{u}\pi,
\ee
with respect to the $L^{2}(g)$ inner product, so that $N$ spans the $g$-normal space to $Im D\pi$. 

  Let $(C_{+}^{m-2,\a})' \subset C_{+}^{m-2,\a}$ denote the subspace of functions $\chi$ such that   
\be \label{norm}
\int_{S^{2}}\chi dv_{+1} = 4\pi.
\ee
Of course for $K = \pi(u)$, $\chi := Ke^{2u} \in (C_{+}^{m-2,\a})'$ by the Gauss-Bonnet theorem. Similarly let 
$(C^{m,\a})' = \{u \in C^{m,\a}: \int u dv_{+1} = 0\}$. 

\begin{lemma}
The map 
\be \label{G}
C: (C^{m,\a})' \to (C_{+}^{m-2,\a})',
\ee
$$C(u) = \mu := Ke^{2u} = \pi(u)e^{2u},$$
is a smooth diffeomorphism.  
\end{lemma} 

{\bf Proof:} We first observe that $C$ is surjective, i.e.~$\mu = Ke^{2u}$ may be arbitrarily 
specified. One needs to show that given any $\mu$ with $\int \mu dv_{+1} = 4\pi$, there is $u$ such that 
$$\D u - 1 = -\mu.$$
This is an elementary consequence of the Fredholm alternative. 

Next we claim that $Ker DC = 0$ everywhere. To see this, suppose $\mu' = 0$. Then $\D u' = 0$, so that 
$u' = const$. Since $u \in (C^{m,\a})'$, $\int u' dv_{+1} = 0$, and hence $u' = 0$. It is easy to see that $C$ is of index 0, 
and hence it follows that $C$ is a local diffeomorphism onto $(C_{+}^{m-2,\a})'$. Finally, to see that $C$ is one-to-one, 
suppose $\mu_{1} = \mu_{2}$. Then $\D u_{1} = \D u_{2}$, which implies $u_{1} =  u_{2} + const$, so again 
$u_{1} = u_{2}$. This proves the result. 

{\endproof}

   We next show that the singular set $\S$ of $\pi$ is a stratified space; this is essentially Theorem 1.4. For the purposes 
of this paper, a subset $X$ of $C^{m,\a}$ will be called a stratified space if $X = \cup S_{i}$, where each $S_{i}$ 
is a smooth manifold of codimension $i$ in $C^{m,\a}$ such that the point-set theoretic closure satisfies  
$$\bar S_{k}\setminus S_{k} \subset \bigcup_{i>k}S_{i}.$$

\begin{theorem} 
The singular set $\S \subset C^{m,\a}$ of $\pi$ is a closed stratified space, with strata $\S_{d}$ consisting of submanifolds of 
codimension $d = dim N \geq 1$. 

  In particular, the regular points of $\pi$ are open and dense in $C^{m,\a}$. 
\end{theorem}

{\bf Proof:} Consider the map 
\be \label{L}
L: (C_{+}^{m-2,\a})' \times (C^{m,\a} \setminus \{0\}) \to C^{m-2,\a},
\ee
$$L(\mu, v) = \D v + 2\mu v.$$
Note that by \eqref{dpi}, $L \circ C = -e^{-2u}D\pi$, so at a critical point of $\pi$, 
$$D^{2}\pi = -e^{2u}D(L \circ C).$$ 

We claim that $0$ is a regular value of $L$, so that the inverse image $S = L^{-1}(0)$ is a smooth properly embedded 
submanifold of $(C_{+}^{m-2,\a})' \times (C^{m,\a}\setminus \{0\})$. To do this, one needs to show the linearization 
is surjective and the kernel splits. 

  The map $L$ is linear in the second component $v$, and $D_{2}L$ maps onto $Im(\D + 2\mu)$ equal to 
the orthogonal complement to the kernel $Ker(\D + 2\mu)$. On the other hand, one has 
\be \label{mu'}
D_{1}L(v) = 2\mu'v,
\ee
where $D_{1}$ is the derivative in the first or $\mu$-direction. We claim that the projection of $\mu' v$ onto the 
kernel $Ker(\D + 2\mu)$ is surjective as $\mu'$ ranges over $T(C_{+}^{m-2,\a})'$. The component of $\mu' v$ 
in the kernel is given by 
\be \label{kerpro}
\int_{S^{2}}\mu' v \o dv_{g_{+1}},
\ee
for some unit vector $\o \in Ker(\D + 2\mu)$, so the claim is that for any $v$, there exists $\mu'$ such that \eqref{kerpro} is non-zero. 
The only condition on $\mu'$ is the linearization of \eqref{norm} which gives
\be \label{linnorm}
\int_{S^{2}}\mu' dv_{g_{+1}} = 0.
\ee
There exists $\mu'$ satisfying \eqref{linnorm} with \eqref{kerpro} non-zero provided $\o v$ is not identically one and $v$ is 
not identically zero. The latter holds by assumption (see \eqref{L}), and the former holds since for instance $\o = 0$ somewhere. 
This shows that $DL$ is surjective. 

  The kernel of $DL$ at $(\mu, v)$ is given by pairs $(\mu', w)$ such that 
\be \label{ker}
\D w + 2\mu w = -2\mu' v.
\ee
This is solvable for $w$ if and only if $\mu' v \perp Ker(\D + 2\mu)$ and the solution is unique modulo $Ker(\D + 2\mu)$. 
Thus $Ker DL$ is a closed subspace, naturally isomorphic to $(Ker(\D + 2\mu)^{\perp}\oplus Ker(\D + 2\mu)$. 
It has a closed complement, given by $\mu'v \in Ker(\D + 2\mu)$ with $w \in Ker(\D + 2\mu)^{\perp}$

  This shows that the set $S = \{(\mu, \o)\} \subset {\rm domain} \, L$ such that $\D \o + 2\mu \o = 0$, 
is a smooth properly embedded submanifold of the domain of $L$. The tangent space $TS$ at $(\mu, w)$ is given by $Ker DL$, 
i.e.~pairs $(\mu', w)$ satisfying \eqref{ker}. Of course $S$ may have a number of distinct components. The reason for excluding 
$\{0\}$ in the domain of $L$ in \eqref{L} is to remove a trivial region $(\mu, 0)$ of $S$, diffeomorphic to $(C_{+}^{m-2,\a})'$. 

   The structure of the singular set $\S$ is now obtained by studying the projection of $S$ into the base space 
$(C^{m,\a})'$. Observe that the vertical fiber of $S$ where $\mu$ is fixed and $\o$ varies is the linear space 
$N$ of dimension $d(\mu)$, so that $S$ is ruled by linear spaces, equal to the fibers $S\cap pr_{1}^{-1}(pt)$ 
of the projection $pr_{1}$ onto the first factor in \eqref{L}. The function 
\be \label{d}
d: S \to \bZ^{+}, \ \ d = dim N(\mu)
\ee
is upper semi-continuous; it may increase in limits, but not decrease. Let $\S = pr_{1}(S)$. 

Let $S_{d} = \{\mu \in S: d(\mu) = d\}$. One should be aware in the following that any given $S_{d}$ may be empty. 
For instance $S_{0} = \emptyset$, since $(\mu, v) \in S_{0}$ implies $v = 0$, which is not in the domain of $L$. 
We proceed inductively on $d$ and work on a fixed component of $S$. First, the upper semi-continuity of $d$ 
implies that $S_{1}$ is open in $S$. Since $S_{1}$ is ruled by lines, the projection
$$pr_{1}: S_{1} \to \S_{1} \subset (C_{+}^{m-2,\a})'$$ 
with $\S_{1} = pr_{1}(S_{1})$ is a vector bundle of rank $1$ (a line bundle with the zero-section removed) over the manifold $\S_{1}$. 
By \eqref{ker}, the space $\S_{1}$ is of codimension 1 in $(C_{+}^{m-2,\a})'$. The upper semi-continuity of $d$ implies 
that $\partial S_{1} = \bar S_{1} \setminus S_{1} \subset \cup_{k\geq 2}S_{k}$.  

  Next, in the closed complement $S\setminus S_{1}$, $S_{2}$ is an open (possibly empty) subset in the induced topology, again 
by the upper semi-continuity of $d$. The projection
$$pr_{1}: S_{2} \to \S_{2} \subset (C_{+}^{m-2,\a})',$$ 
with $\S_{2} = pr_{1}(S_{2})$ is a vector bundle of rank $2$ over the manifold $\S_{2}$. Again by \eqref{ker}, the manifold 
$\S_{2}$ is of codimension 2 in $(C_{+}^{m-2,\a})'$. 

  Continuing inductively in the same way gives manifolds $\S_{d}$ of codimension $d$ in $(C_{+}^{m-2,\a})'$ such that  
the point-set theoretic boundary of each $\S_{d}$ is contained in $\cup_{d'>d}\S_{d'}$. This gives $\S$ the structure 
of a stratified space. 

  Via the diffeomorphism $C$ from Lemma 5.4, we view $\S = pr_{1}(S)$ as a subset of $\cU'$. By scaling 
(i.e.~ the addition of constants to $u$), $\S$ becomes a stratified subset of $\cU$ with strata $\S_{d}$ of 
codimension $d$ in $C^{m,\a}$. 

{\endproof}

  Next we turn to discussion of the singular value set $\cZ \subset C_{+}^{m-2,\a}$. By definition, 
$$\cZ = \pi(\S).$$  
The singular set $\S$ in the domain $C^{m,\a}$ is a stratified space, with strata $\S_{d}$ given by submanifolds of codimension $d$ 
in $C^{m,\a}$ with (point-set theoretic) boundary satisfying $\partial \S_{d} \subset \S_{d+1}$. The map $\pi$ maps $\S$ onto 
the set of singular values $\cZ$ in $C_{+}^{m-2,\a}$ so that $\cZ$ is a singular ``variety" of codimension at least one in 
$C_{+}^{m-2,\a}$. More precisely, $\cZ$ is a rectifiable set in $C_{+}^{m-2,\a}$ of codimension at least one, (cf.~the definition 
following Proposition 2.13). 

   Let 
$$\cZ_{d} = \pi(\S_{d}),$$
so that, if non-empty, $\cZ_{d}$ is rectifiable of codimension $d$ in $C_{+}^{m-2,\a}$. It is not known if $\cZ$ is stratified as is 
$\S$. The following result gives a sufficient criterion for the strata $\cZ_{d}$ to be locally given by submanifolds of codimension $d$ in 
$C_{+}^{m-2,\a}$. 

\begin{lemma}
Suppose $u \in \S_{d}$ with $K = \pi(u) \in \cZ_{d}$ and let $N$ be the $d$-dimensional kernel of $\D + 2Ke^{2u}$. 
If for some $\w \o \in N$, the symmetric bilinear form $S_{\w \o}: N\times N \to \bR$, 
\be \label{N}
S_{\w \o}(\o_{1}, \o_{2}) = \int_{S^{2}}K\o_{1} \o_{2}\w \o dv_{g},
\ee
is non-degenerate, then $\cZ_{d}$ is a submanifold of codimension $d$ near $K$. 
\end{lemma}

{\bf Proof:} 
 For $u \in \S_{d}$, consider the straight line $u_{t} = u + tu'$, where $u' = \o \in N$, so $\o$ is normal to $Im D_{u}\pi$ in 
$L^{2}(g)$. Then $\s_{t} = \pi(u_{t}) = K(u_{t})$ is a smooth curve in $C_{+}^{m-2,\a}$, with vanishing tangent vector at $t = 0$, 
$\s'(0) = D\pi(u') = 0$. Calculating the second derivative from \eqref{dpi} gives, since $u'' = 0$,  
\be \label{two}
D^{2}\pi(\o, \o) = 2e^{-2u}\D\o = -4K\o^{2}. 
\ee
Thus taking the $L^{2}(g)$ pairing with the normal vector $\w \o$ as 
in \eqref{perp} gives 
$$\<D^{2}\pi(\o, \o), \w \o\> = -4\int_{S^{2}}K\o^{2} \w \o dv_{g}.$$

  The result is then a straightforward consequence of the implicit function theorem, cf.~Theorem 3.7.2 of \cite{N2} 
for instance. 

{\endproof} 
 
  Next we consider relations between the singular set $\cZ$ and the existence of solutions of \eqref{u}. As discussed in 
Sections 2 and 4, the boundary $\partial \cK$ is a closed rectifiable set of codimension 1, separating $C_{+}^{m-2,\a}$ into 
an infinite collection of domains partially distinguished by the degree of $\pi_{0}$. There are only three possible regions 
where \eqref{u} may not have solutions. The first region is the degree zero component $\cK^{0}$, cf.~Remark 4.2. By Proposition 
4.7, (i.e.~Theorem 1.5), a second region is a portion of the boundary region $\partial \cC \subset \partial \cK$, while the 
third region is the boundary $\partial \cN \subset \partial \cK$. We explore existence and non-existence in these 
regions in more detail. 

  We begin first with $\cK^{0}$. For $\pi^{0}: \cU^{0} \to \cK^{0}$ as in \eqref{pii2}, let 
$$\cE = Im \pi^{0} \subset \cK^{0}.$$
\begin{proposition}
The image $\cE$ of $\pi^{0}$ is a non-empty, closed subset of $\cK^{0}$. 
\end{proposition}

{\bf Proof:} It follows from Proposition 4.3 that $\cE \neq \emptyset$. Since $\pi^{0}$ is proper, it has closed range, 
which proves the second statement. 

{\endproof}

 Let $\cE_{int}$ denote the interior of $\cE$, as a set in $C_{+}^{m-2,\a}$; this is the union of the set of regular values of 
$\pi^{0}$ with the singular values (if any) contained in the interior, i.e.~ ($\cE_{int}\cap \cZ$). Obviously $\cE_{int}$ is open 
in $\cE$. More significantly, it follows from Theorem 5.5 and the fact that $\pi^{0}$ is smooth that $\cE_{int}$ is dense in $\cE$, 
\be \label{kdense}
\bar \cE_{int} = \cE.
\ee

   Let $\cB = \partial \cE_{int}$ be the boundary of $\cE_{int}$, so that $\cE = \cE_{int} \cup \cB$. One 
thus has a decomposition of the target 
\be \label{decomp}
\cK^{0} = \cE_{int}\cup \cB \cup \cD,
\ee
where $\cD = \cK^{0} \setminus \cE$ is the non-existence region where there are no solutions to the Nirenberg problem. 
The linear functions $K = \ell$ are in $\cK^{0}$ but are not in $Im \, \pi^{0}$ by the Kazdan-Warner obstruction \eqref{kw} 
so that 
$$\cD \neq \emptyset.$$
Both $\cE_{int}$ and $\cD$ are non-empty open sets, so that the decomposition \eqref{decomp} is non-trivial. 
Of course the boundary satisfies 
$$\cB \subset \cZ.$$

\begin{proposition} 
The wall $\cB$ separating $\cE_{int}$ from $\cD$ satisfies 
\be \label{codim1}
\cB \subset \bar \cZ_{1}.
\ee
In particular, $\cZ_{1} \neq \emptyset$. An open-dense subset $\cB'$ of $\cB_{1} = \cB\cap \cZ_{1}$ is a smooth 
hypersurface in $\cK^{0} \subset C_{+}^{m-2,\a}$. 
\end{proposition}

{\bf Proof:} The singular strata $\S_{k}$ are of codimension $k$ in the domain $\cU^{0}$. Hence their images $\cZ_{k}$ under 
the smooth map $\pi$ are of codimension at least $k$ in $C_{+}^{m-2,\a}$. Closed subsets of codimension at least two 
do not locally separate the space $C_{+}^{m-2,\a}$, and since $\cB$ does locally separate, one must have 
$\cZ_{1} \neq \emptyset$ and \eqref{codim1} follows. 

   We use Lemma 5.6 to determine which regions in $\cB_{1}$ are codimension 1 submanifolds. First, define a point $K \in \cB_{1}$ to be 
{\it  path-accessible} if $\cE_{int}$ is locally path connected near $K$. These are points $K$ such that for any open neighborhood $V$ 
of $K$ in the target $C_{+}^{m-2,\a}$, there is a neighbhorhood $U \subset V$ of $K$ such that $U \cap \cE_{int}$ is path connected. 
This fails for instance for points $K$ at which two local regions of $\cB_{1}$ intersect at $K$, separating $\cE_{int}$ into two 
local components near $K$. Let $\cB'$ denote the set of path-accessible points in $\cB_{1}$. 

  By construction, $\cB'$ is open in $\cB_{1}$. We claim that $\cB'$ is also dense in $\cB_{1}$. Namely for any point 
$K \in \cB_{1}$ and for any neighborhood $U$ of $K$ in $C_{+}^{m-2,\a}$, one has $U \cap \cD \neq \emptyset$. It follows 
from this that $U\cap \cB' \neq \emptyset$, so that $\cB'$ is dense. 

  Next we prove that the path-accessible points $\cB'$ form (locally) a smooth hypersurface in $C_{+}^{m-1,\a}$. To do this, it 
suffices to show that the form $S$ in \eqref{N} is non-degenerate on $\cB'$. 

   Suppose the form $S = D^{2}\pi$ is degenerate at $u$ with $\pi(u) = K \in \cB'$. Then there exists $u' = \o \in N$ 
and $\w \o \in N$ such that $\<D^{2}\pi(\o, \o), \w \o\> = 0$. We then calculate the third derivative $D^{3}\pi(\o, \o, \o) = 
\frac{d^{3}(K_{u_{t}})}{dt^{3}}$. Taking the derivative of \eqref{two} along $u_{t}$, using the fact that $K'(\o) = 0$ and $u'' = 0$ 
gives 
$$D^{3}\pi(\o, \o, \o) = -4e^{-2u}\o^{2}\D \o = 8K\o^{3},$$
so that 
\be \label{three}
\<D^{3}\pi(\o, \o, \o), \w \o\> =  8\int_{S^{2}}K\o^{3} \w \o dv_{g}.
 \ee
When $dim N = 1$, $\w \o$ is a multiple of $\o$. Using \eqref{dpi} (with $v = \w \o = \o$) gives 
$$\<D^{3}\pi(\o, \o, \o), \o\> = -4\int_{S^{2}}\o^{3}\D_{g}\o dv_{g} = 12\int_{S^{2}}\o^{2}|d\o|^{2} dv_{g} > 0.$$

  Now as in the proof of Lemma 5.6,  suppose $\s(0) = \pi(u_{0})= K \in \cB_{1}$, so that $u \in \S_{1}$. The normal space to $Im D\pi$ at $K$ 
is 1-dimensional and spanned by $\o$. The term \eqref{three} is then strictly positive and hence the third derivative of $\s(t)$ in the direction 
of the normal vector $\o$ is non-vanishing. It follows that the function $\<\s(t), \o\>$ behaves as a cubic, and so changes signs. 
Moreover, by construction $\s(t) \in \cE$ and since $\s'(t) \neq 0$ for $t \neq 0$, it follows that $\s(t) \in \cE_{int}$ for 
$t \neq 0$ with $t$ small. Now if $K$ is path-accessible, then $\s(t)$ is in a single local path component of $\cE_{int}$, which 
contradicts the fact that $\s(t)$ behaves as a cubic. This contradiction shows that the form \eqref{N} is non-degenerate on $\cB'$, 
which proves the result.

{\endproof}

\begin{remark}
{\rm The open, dense set $\cB' \subset \cB$ where $\cB$ is a submanifold of codimension 1 in $C_{+}^{m-2,\a}$ is a classical 
{\it bifurcation locus} for the Nirenberg problem, cf.~\cite{N2} for instance. Namely for $K \in \cB'$ with (outward) normal vector $\o$, 
consider the line $K(t) = K + t\o$ in $C_{+}^{m-2,\a}$. For $t > 0$ sufficiently small, $K(t) \notin Im \,\pi$, i.e.~there are no solutions 
of \eqref{u} for such $K(t)$. For $t < 0$ sufficiently small, there are locally exactly two distinct solutions of \eqref{u}, i.e.~the map $\pi$ 
is locally $2-1$ onto $K(t)$, $t < 0$. These two solutions merge at $t = 0$ to give a local unique solution to the Nirenberg problem. 
Thus the map $\pi$ is locally a $2-1$ fold map onto the curve $K(t)$, $t \leq 0$. 

  A large open set in the non-existence region $\cD$ is given by the curvature functions from Proposition 4.4. 
A further specific example of non-existence of solutions in $ \cD \subset \cK^{0}$ near the boundary with the region 
$\cK^{1}$ (a degree one component) is constructed and analysed in detail by Struwe in \cite{St}. 

}
\end{remark}

\begin{remark}
{\rm As an application of the results above, we discuss the behavior of the map $\pi^{0}$ near the special point 
$K = 1$ (or $K = const$). Consider the region $\cP$ of positive curvature functions in $C_{+}^{m-2,\a}$, as in \eqref{cP}. 
Then $\cP \subset \cC$ and as noted in \cite{CGY}, the domains $\cK^{i}\cap \cP$ are star-shaped about $K = 1$, i.e.~if $K 
\in \cK^{i}\cap \cP$, then the line segment $t+ (1-t)K$, $t \in [0,1)$ lies in $\cK^{i}\cap \cP$. Thus the structure 
of $\partial \cK$ becomes very complicated near $K = 1$. 

   In the component $\cK^{0}$ where  
$$deg \pi^{0} = 0,$$
any regular value $K \in \cE_{int}$ near $K = 1$ has an even number of solutions of \eqref{u}. On the other hand, 
it is easy to see that the point $u = 0$ with $\pi(u) = K_{u} = 1$ is a regular point for the restricted map 
$\pi_{1} = \pi|_{\cS}$ for $\cS$ as in \eqref{cS}, i.e.~$D_{0}\pi_{1}$ is injective. Points $u \in \cS$ near $u = 0$ 
are then also regular and so in a neighborhood $V$ of $0$ in $\cS$, the map 
\be \label{pi1}
\pi_{1}: V \to H,
\ee
is a diffeomorphism onto its image $H$; $H$ is a local smooth hypersurface of codimension 3 in $C_{+}^{m-2,\a}$. 

  Consider then a regular value $K \in H\cap \cK^{0}$. The value $K$ has a unique inverse image (i.e.~solution of \eqref{u}) 
$u_{1} \in V \subset \cS$ and so there must be another solution $u_{2} \notin V$ of \eqref{u}. (In fact with further 
work it can be proved that $u_{2} \notin \cS$). Moreover, the solutions $u_{1}$ and $u_{2}$ have $J_{K}$-indices 
with opposite (even-odd) parity. By \eqref{cslice}, there exists a conformal dilation $\f$ such that 
$$\f^{*}u_{2} \in \cS,$$ 
and so 
$$\pi(\f^{*}u_{2}) = K\circ \f.$$ 
Suppose that $\w u_{2} = \f^{*}u_{2} \in V$, so that $K \circ \f \in H$. Then $u_{1}$ and $\w u_{2}$ are two regular 
points of $\pi^{0}$ in $V \subset \cS$ with index of opposite parity. However, this is not possible. Namely, one may 
take a curve $K(t)$ from $K$ to $K\circ \f$ in $H$ with $K(t)$ a regular value of $\pi$ for all $t$. The curve $K(t)$ 
has a unique lift to a curve $u(t) \subset V$. As discussed in the proof of Theorem 1.3 in Section 3, the index of 
$u(t)$ is then independent of $t$, giving a contradiction. It follows that 
\be \label{nV}
\w u_{2} \notin V.
\ee

  This argument proves that if $K \in \cS$ is sufficiently close to $1$, i.e.~$|K -1|_{C^{m-2,\a}} \leq \e$ for $\e$ sufficiently 
small, then 
$$|\f|_{C^{1}} \geq C = C(\e),$$
where $C(\e) \to \infty$ as $\e \to 0$. For if $C$ were uniformly bounded as $\e \to 0$, then the solution $\w u_{2}$ above 
would necessarily be in $V$, contradicting \eqref{nV}. 

   Thus, pairs of solutions to \eqref{u} with $K \in \cK^{0}$ sufficiently close to 1 are (very) far apart in $C^{m,\a}$, 
differing from each other by (very) large conformal factors. This illustrates concretely the effect of the non-compactness 
of the conformal group. 

}
\end{remark}

   Proposition 5.8 describes the general structure of the wall separating the regions of existence and non-existence in 
$\cK^{0}$. We next consider briefly the issue of existence and non-existence near the second region $\partial \cK$. 
Recall that 
$$\partial \cK = \partial \cC \cup \partial \cN.$$
By Proposition 4.7, for any $(m, n)$ with $m < n$ and $m \neq 0$, one has
$$\partial \cC_{m}^{n} \subset Im \, \pi,$$
so that for any $K \in \partial \cC_{m}^{n}$ there is a solution $u$ of \eqref{u} with $K_{u} = K$. From the proof of 
Proposition 4.7, for any such $(m, n)$, say $n = m+1$, there are (large families of) solution curves $u(t)$, $t \in [0,1]$ 
with $K_{u(t)} \in \cK^{m}$ for $t < \frac{1}{2}$, $u(\frac{1}{2}) \in \partial C_{m}^{n}$ and $K_{u(t)} \in \cK^{n}$ for 
$t > \frac{1}{2}$. Thus, curves of solutions do not necessarily ``blow-up" on approach to $\partial \cC$ and $\pi$ 
is ``proper" along such curves. 

  We expect similar behavior holds at $\partial \cN$ but this remains unknown in general. The eigenfunction examples 
discussed in Proposition 4.5 and Remark 4.6 show that 
$$Im \, \pi \cap \partial \cN \neq \emptyset.$$
Moreover, Remark 5.10 holds on any component $\cK^{n}$ of degree $n$ and not just on $\cK^{0}$; for any such $\cK^{n}$, 
there are solutions $u$ with $u \in \cK^{n}\cap H$, for $H$ as in \eqref{pi1}. This follows for instance from the work of 
Chang-Gursky-Yang \cite{CGY}, where $deg \, \pi^{n}$ is shown to be the intersection number of the orbit $\{K\circ \f\}$, 
$\f \in \Co$ with $H$, for $K$ sufficiently close to $1$. As above with $\partial \cC$, there are (large families of) solution curves 
$u(t) \subset V \subset \cS$ with $K_{u(t)}$ passing smoothly from one component of $\cK^{m}\cap \cP$ to another 
component $\cK^{m+1}\cap \cP$ within $H$.

\begin{remark}
{\rm Since the degree of $\pi_{0}$ changes on passing from one component of $\partial \cK$ to a neighboring or bordering 
component, $\pi_{0}$ cannot extend to a proper map on any open set containing a portion of $\partial \cK$. In particular, 
the {\it apriori} estimates \eqref{apri} do not extend to such larger domains. On curves $K(t)$ passing through $\partial \cC 
\subset \partial \cK$, there must exist curves of solutions $u(t)$ which blow-up as $K_{u(t)}$ approaches $\partial \cK$. 
The analysis given in the proof of Proposition 4.7 shows that in a passage from $\cK^{m+1}$ to $\cK^{m}$, $m \neq 0$, 
it is likely that only one solution curve blows up, while the remaining $m$ solution curves pass continuously through 
$\partial \cC_{m}^{m+1}$, i.e.~a single blow-up curve is the reason for the drop in degree from $m+1$ to $m$. Similarly, 
in passing from $\cK^{1}$ to $\cK^{0}$ through $\partial \cC$ one expects generically there is only one solution curve,
and such a solution curve blows up on approach to $\partial \cC_{0}^{1}$. 

One could expect similar behavior on passing through the $\partial \cN$ region of $\partial \cK$.

}
\end{remark}

\section{Symmetry Breaking and Existence}

   In this section, we turn to the question of existence of solutions of \eqref{u} with extra symmetry. A key point of 
view is to find conditions which break the symmetry of the (non-proper part of the) action of the conformal group 
${\rm Conf}(S^{2})$ on the domain and target spaces of the map $\pi$. 

\medskip 

  To begin, we give a new proof of Moser's theorem \cite{M2}, that any even function $K$ on 
$S^{2}$ ($K(-x) = K(x)$) with $K > 0$ somewhere is the Gauss curvature of a conformal metric on $S^{2}$. 
The proof below does not use any sharp Sobolev inequality of Moser-Trudinger-Aubin type. 

  Let $A$ be the antipodal map $x \to -x$ of $S^{2}$. Let $C_{ev}^{m,\a} \subset C^{m,\a}$ be the subspace of 
even functions, so $u(-x) = u(x)$ and similarly consider $C_{ev}^{m-2,\a} \subset C_{+}^{m-2,\a}$. The restriction of 
$\pi$ in \eqref{pi} to $C_{ev}^{m,\a}$ gives a smooth Fredholm map 
\be \label{piev}
\pi_{ev}: C_{ev}^{m,\a} \to C_{ev}^{m-2,\a}.
\ee
Let $\cC_{ev} = \cC \cap C_{ev}^{m-2,\a}$. 

\begin{theorem}  
Any even function $K \in C_{+}^{m-2,\a}$ is the Gauss curvature of a metric $\g = e^{2u}\g_{+1}$, 
with $u$ an even function on $S^{2}$. Moreover, the restriction of $\pi_{ev}$ in \eqref{piev} to $\cU_{ev} = 
\pi_{ev}^{-1}(\cC_{ev})$ is proper and the degree (up to sign) of the map 
\be \label{pievc}
\pi_{ev}: \cU_{ev} \to \cC_{ev}, 
\ee
satisfies  
\be \label{degev}
deg \, \pi_{ev} = 1.
\ee
\end{theorem}

{\bf Proof:} As before, via \eqref{u}, it is easy to see that $\pi_{ev}$ in \eqref{piev} is Fredholm of index 0. 
Recall from the work in Section 2 that the non-properness of the map $\pi$ has two sources. First, the 
non-compactness of the conformal group $\Co$,  which gives rise to the division $\partial \cN$ into the regions 
$\cN$ and second the `bubble formation' at $\partial \cC$, giving rise to the region $\cC$. We will deal with 
these in turn. 

   To begin, the conformal group of $\bR \bP^{2}$ is {\it compact}, equal to its isometry group. Propositions 2.1, 2.8 and 
Theorem 2.3 apply equally well to $\bR \bP^{2}$ in place of $S^{2}$ and the proof of Theorem 2.9 for the map $\pi_{ev}$ in 
\eqref{pievc} becomes much simpler. The compactness of the conformal group implies that the diffeomorphisms 
$\psi_{i}$ of $\bR \bP^{2}$ in \eqref{a1}, (which lift to even diffeomorphisms of $S^{2}$) themselves converge (in a 
subsequence) to a limit $\psi$, so that the renormalized diffeomorphisms $\psi_{i}$ in \eqref{a3} converge to the 
identity in $C^{m+1,\a}$. This shows that $\pi_{ev}$ in \eqref{pievc} is indeed proper. In effect, $\partial \cN = 
\emptyset$ and $\cN = C_{ev}^{m-2,\a}$ in this situation. 

   Proposition 2.14 shows that $\cC$ is connected and the same proof shows that $\cC_{ev}$ is connected. Hence 
$deg \, \pi_{ev}$ is well-defined on $\cU_{ev}$. 

  To compute the degree of $\pi_{ev}$, consider the value $K = 1$. The solutions of \eqref{u} with $K = 1$ are given by 
$u = \log \chi$, for $\chi$ as in \eqref{confact}. The only such $u$ which is even, and so descends to $u$ on $\bR \bP^{2}$ is 
$u = 0$. Thus, $\pi_{ev}^{-1}(1)$ is the function $u = 0$. Moreover, the kernel $Ker D\pi$ of $\pi$ at $u = 0$ is given 
by the linear functions $\ell$, of which the only even function is $\ell = 0$. Thus $Ker D\pi_{ev} = 0$. This shows that 
$K = 1$ is a regular value of $\pi_{ev}$, uniquely realized by the regular point $u = 0$. This proves \eqref{degev}. 

  It follows that $\pi_{ev}$ is onto $\cC_{ev}$. Theorem 1.5, (or more precisely its proof in Proposition 4.7), shows 
that $\pi_{ev}$ is also onto $\partial \cC_{ev}$. (There are no components of $\cK = \cN \cap \cC = \cD$ of degree 
zero. This completes the proof. 

{\endproof}

\begin{remark}
{\rm It is not asserted (and is not true) that $\pi_{ev}$ in \eqref{piev} is proper. The formation of bubbles as discussed 
in Remark 2.7 holds on $\bR \bP^{2}$ just as well as on $S^{2}$. It follows in particular that the number of even solutions 
of \eqref{u} as $K \in \cC_{ev}$ approaches a given function in $\partial \cC$ from the degree decreasing side (so a local 
maximum of $K$ transitions from positive to negative values) is often different than the number of even solutions as $K$ 
approaches a given function from the degree increasing side of $\partial \cC$; compare with the proof of Proposition 4.7.  

}
\end{remark}

 It is not difficult to see that the result above generalizes, with essentially the same proof, to the action of any compact 
group $\Gamma \subset O(3) = Isom(S^{2})$ of isometries of $S^{2}$ which breaks the non-compact action 
of $\Co$ on the target space $C_{+}^{m-2,\a}$. Thus, we assume $\G \subset O(3)$ satisfies the following: 
for any conformal dilation $\f = \f_{p,t}$, $\f \neq Id$ in $\Co$, there exists $\g \in \G$ such that 
\be \label{nc1}
\g \circ \f \neq \f \circ \g,
\ee
as elements in $\Co$. For example the subgroup $\bZ_{2} \subset O(3)$ generated by the antipodal map 
satisfies \eqref{nc1}. 

  Let $C_{\G}^{m-2,\a}$ be the space of $\G$-invariant functions in $C_{+}^{m-2,\a}$, i.e.~$K(\g(x)) = K(x)$, 
for all $x \in S^{2}$, $\g \in \G$ and similarly for $C_{\G}^{m,\a}$. As before, $\pi$ restricts to a smooth Fredholm 
map of index 0, 
$$\pi_{\G}: C_{\G}^{m,\a} \to C_{\G}^{m-2,\a}.$$
As in \eqref{pievc}, consider also 
\be \label{piG}
\pi_{\G}: \cU_{\G} \to \cC_{\G},
\ee

  We then have: 

\begin{theorem} 
If a finite group action $\G \subset O(3)$ as above breaks the non-compact action of the conformal group, 
then any $\G$-invariant function $K \in C_{+}^{m-2,\a}$ is the Gauss curvature of a conformal metric 
$\g = e^{2u}\g_{+1}$ with $u$ also $\G$-invariant. Further, the map $\pi_{\G}$ in \eqref{piG} is proper 
and (up to sign)  
\be \label{degg}
deg \, \pi_{\G} = 1.
\ee
\end{theorem} 

{\bf Proof:} The proof is identical to the proof of Theorem 6.1, with one distinction. Let $S = S^{2}/\G$ be the orbit space of 
the $\G$-action on $S^{2}$. Then $S$ is an orbifold, in fact a ``good" orbifold in the sense of Thurston. A classification of 
orbifolds in the special case when $\G \subset SO(3)$ is given for instance in \cite{Z}.  As in the proof of Theorem 6.1, the metrics 
$g_{i}$ descend to orbifold metrics $\w g_{i}$ on $S$, and one needs an orbifold or $\G$-equivariant version of 
Propositions 2.1, 2.8  and Theorem 2.3 to complete the proof as before. 

   However, since the $\G$-action is fixed for all metrics $g_{i}$, this is essentially standard. Namely, recall that the proof 
of the Gromov convergence theorem (Proposition 2.1) is fundamentally based on the construction of suitable harmonic coordinate 
charts of uniform size in which the metric is well-controlled, cf.~\cite{GW} or \cite{Pe}. Since $\G$ acts by isometries on 
$(S^{2}, g_{i})$, for each $i$ one may choose a suitable $\G$-equivariant finite atlas $\cA$ of harmonic coordinates on $(S^{2}, g_{i})$, 
i.e.~$\{x_{i}^{k}\}: U \subset S^{2} \to B \subset \bR^{2}$ is a $g_{i}$-harmonic coordinate chart in $\cA$ if and only if 
$\{x_{i}^{k}\circ \g\}$ is in $\cA$, for each $\g \in \G$. The collection of transition maps between charts in $\cA$ are thus also 
$\G$-equivariant. The geometric center of mass construction used in \cite{GW} or \cite{Pe} to pass from local to global maps is also 
$\G$-equivariant. Given this, it is readily verified by inspection that the proof in \cite{GW} or \cite{Pe} implies the 
$\G$-equivariant convergence of the metrics modulo diffeomorphism, i.e.~the diffeomorphisms $\psi_{i}$ as in \eqref{a1} 
satisfy 
$$\g \circ \psi_{i} = \psi_{i}\circ \g,$$ 
for all $\g \in \G$ and so are orbifold diffeomorphisms. The condition \eqref{nc1} ensures that the group of $\G$-equivariant 
conformal diffeomorphisms is trivial, so that.by the same arguments used before in the proof of Theorem 6.1, $\pi_{\G}$ in 
\eqref{piG} is proper. The remainder of the proof then proceeds exactly as in the proof of Theorem 6.1. 

{\endproof}

 As a simple example, one may take $\Gamma$ to be the group generated by rotation by angle $\pi$ (or any other non-zero 
rational angle) along one axis and rotation by angle $\pi$ (or any other non-zero rational angle) along a linearly independent 
axis.

\begin{remark}
{\rm The surjectivity statements in Theorems 6.1 and 6.3 also hold for thickenings of $C_{\G}^{m-2,\a}$, i.e.~open sets in 
$C_{+}^{m-2,\a}$ containing $C_{\G}^{m-2,\a}$. Namely, first observe that $C_{\G}^{m-2,\a}$ is not contained in the singular 
value set $\cZ$. To see this, suppose $K \in C_{\G}^{m-2,\a} \cap \cZ$ and let $\o$ be an outward normal vector to $T\cZ$ 
in $C_{+}^{m-2,\a}$, i.e.~$\o \in N$; $\o$ is generically unique up to scalar multiples. Choose a $\G$-invariant function 
$\chi$ such that 
$$\<\chi, \o\> = \int_{S^{2}}\chi \o dv_{+1} \neq 0.$$
Then there exists $\e$ small such that $K_{\e} = K - \e \chi \in C_{\G}^{m-2,\a} \setminus \cZ$. 

   Thus, although $C_{\G}^{m-2,\a}$ is of infinite codimension in $C_{+}^{m-2,\a}$, there are open neighborhoods $\cV$ of 
$C_{\G}^{m-2,\a}$ contained in $Im \, \pi$. Hence the full ${\rm Conf}(S^{2})$ orbit of $\cV$ is contained 
in $Im \, \pi$. 

}
\end{remark} 

\begin{remark}
{\rm Although $deg \, \pi_{ev}$ or $deg \, \pi_{\G}$ in \eqref{degev} or \eqref{degg} equals one (up to sign), this does not 
imply that the image $C_{\G}^{m-2,\a}$ is contained in the degree one component $\cN^{1}$ (or $\cN^{-1}$) of $\pi_{0}$. 
For instance, for $\G = \bZ_{2}$ it is easy to see from the degree formula \eqref{degform0} that for any odd $n \in \bZ$, 
there exist $K \in C_{ev}^{m-2,\a}$ such that $K \in \cN^{n}$. Hence, for such $K$, there are at least $n$ distinct solutions 
of \eqref{u}; one expects that often only one these $n$ solutions is an even function. 

   Finally, Remark 6.2 also applies to the map $\pi_{\G}$. 

}
\end{remark}

 \bibliographystyle{plain}

\end{document}